\input ./style/arxiv-general.cfg
\documentclass[MSNbibl,number,citesort,nameyear,dvips]{arxbj}
\makeatletter
   \@ifpackageloaded{graphicx}{}{\usepackage{graphicx}}
\makeatother

\volume{22}
\issue{3}
\pubyear{2016}
\firstpage{1808}
\lastpage{1838}
\doi{10.3150/15-BEJ712}
\docsubty{FLA}

\makeatletter
\newcommand{\angler}{\rangle}
\newcommand{\anglel}{\langle}
\newcommand{\rrvert}{\vert}
\newcommand{\rrVert}{\Vert}
\newcommand{\llvert}{\vert}
\newcommand{\llVert}{\Vert}
\newcommand{\nnorm}[1]{\vert\!\vert\!\vert #1\vert\!\vert\!\vert}
\newtheorem{theo}{Theorem}[section]
\newtheorem{lemma}[theo]{Lemma}

\newtheorem{prop}[theo]{Proposition}
\newproclaim{defi}[theo]{Definition}
\newremark{algo}{Algorithm}[section]
\newremark{prob}{Open Problem}[section]
\newremark{rem}{Remark}
\newremark{example}{Example}
\makeatother

\begin{document}
\begin{frontmatter}

\title{Markov Chain Monte Carlo confidence~intervals}
\runtitle{Markov Chain Monte Carlo confidence intervals}

\begin{aug}
\author[A]{\inits{Y.F.}\fnms{Yves F.}~\snm{Atchad\'e}\corref{}\ead[label=e1]{yvesa@umich.edu}}
\address[A]{University of Michigan, 1085 South University, Ann Arbor,
48109, MI, United States.\\ \printead{e1}}
\end{aug}

%
\received{\smonth{2} \syear{2014}}
%
\revised{\smonth{10} \syear{2014}}

\begin{abstract}
For a reversible and ergodic Markov chain $\{X_n, n\geq0\}$ with
invariant distribution $\pi$, we show that a valid confidence interval
for $\pi(h)$ can be constructed whenever the asymptotic variance
$\sigma^2_P(h)$ is finite and positive. We do not impose any additional
condition on the convergence rate of the Markov chain. The confidence
interval is derived using the so-called fixed-b lag-window estimator of
$\sigma_P^2(h)$. We also derive a result that suggests that the
proposed confidence interval procedure converges faster than classical
confidence interval procedures based on the Gaussian distribution and
standard central limit theorems for Markov chains.
\end{abstract}

\begin{keyword}
\kwd{Berry--Esseen bounds}
\kwd{confidence interval}
\kwd{lag-window estimators}
\kwd{martingale approximation}
\kwd{MCMC}
\kwd{reversible Markov chains}
\end{keyword}
\end{frontmatter}

\section{Introduction}\label{secintro}

Confidence intervals play an important role in Monte Carlo simulation
(Robert and Casella \cite{robertetcasella04}, Asmussen and Glynn \cite{asmussen07}). In Markov Chain Monte Carlo
(MCMC), the existing literature requires the Markov chain to be
geometrically ergodic for the validity of confidence interval
procedures (Jones \textit{et~al.}  \cite{jonesetal06},   Flegal and Jones \cite{jonesetal09},
Atchad{\'e} \cite{atchadeaos11}). The main
objective of this work is to simplify some of these assumptions. We
show that for a reversible ergodic Markov chain, a~valid confidence
interval can be constructed whenever the asymptotic variance itself is
finite. No additional convergence rate assumption on the Markov chain
is required.

Let $\{X_n, n\geq0\}$ be a reversible stationary Markov chain with
invariant distribution $\pi$. For $h\in L^2(\pi)$, the asymptotic
variance of $h$ is denoted $\sigma_P^2(h)$ (see (\ref{sig}) below for
the definition). A remarkable result by C. Kipnis\vspace*{1pt} and S. R. Varadhan
(Kipnis and Varadhan \cite{kv86}) says that if $0<\sigma_P^2(h)<\infty$, then $\frac
{1}{\sigma_P(h)\sqrt{n}}\sum_{i=1}^n (h(X_i)-\pi(h))$ converges weakly\vspace*{-2pt}
to $\mathbf{N}(0,1)$ where $\pi(h)\stackrel{\mathrm{def}}{=}\int
h(z)\pi(\mathrm{d}z)$. In order
to turn this result into a confidence interval for $\pi(h)$, an
estimator $\sigma_n$ of $\sigma_P(h)$ is needed. A common practice
consists in choosing $\sigma_n$ as a consistent \mbox{estimator} of $\sigma
_P(h)$. However, consistent estimation of $\sigma_P(h)$ typically
requires further assumptions on the convergence rate of the Markov
chain (typically geometric ergodicity), and on the function~$h$.
Instead of insisting on consistency, we consider the so-called fixed-b
approach developed by Kiefer, Vogelsang and Bunzel \cite{kieferetal00},  Kiefer and Vogelsang \cite{kieferetvogelsang05}, where
the proposed estimator $\sigma_n$ is known to be inconsistent. Using
this inconsistent estimator we show in Theorem~\ref{thm1} that a
Studentized analog of the Kipnis--Varadhan's theorem holds: if
$0<\sigma
_P^2(h)<\infty$, then $\mathbf{T}_n\stackrel{\mathrm{def}}{=}\frac
{1}{\sigma_n\sqrt
{n}}\sum_{i=1}^n (h(X_i)-\pi(h))$ converges weakly to a (non-Gaussian)
distribution. The theorem extends to nonstationary Markov chains that
satisfy a very mild ergodicity assumption. To a certain extent, the
result is a generalization of Atchad{\'e} and Cattaneo \cite{atchadeetcattaneo14} which
establishes the same limit theorem for geometrically ergodic (but not
necessarily reversible) Markov chains. The result is particularly
relevant for Markov chains with sub-geometric convergence rates. For
such Markov chains, the author is not aware of any result that
guarantees the asymptotic validity of confidence intervals. However, it
is important to point out that the finiteness of $\sigma^2_P(h)$
carries some implications in terms of convergence rate of $P$, and is
not always easy to check. But the main point of this work is that the
finiteness of $\sigma^2_P(h)$ is all that is needed for consistent
confidence interval.

As we shall see, Theorem~\ref{thm1} comes from the fact that there
exists a pair of random variables $(N,D)$, say, such that the joint\vspace*{-2pt}
process $ (\frac{1}{\sqrt{n}}\sum_{i=1}^n(h(X_i)-\pi(h)),
\sigma
_n^2 )$ converges weakly to $(\sigma_P(h)N,\sigma_P^2(h)D)$. As a
result, $\sigma_P(h)$ cancels out in the limiting distribution of
$\mathbf{T}_n$. This approach to confidence intervals is closely
related to the standardized time series method of Schruben \cite{schruben83}
(see also Glynn and Iglehart \cite{glynnetiglehart90}), well known in operations research.
Indeed in its simplest form, the standardized time series method is the
analog of the fixed-b procedure using the batch-mean estimator with a
fixed number of batches. Despite this close connection, this paper
focuses only on the fixed-b confidence interval.

We also compare the fixed-b lag-window estimators with the more
commonly used lag-window estimators. We limit this comparison to the
case of geometrically ergodic Markov chains. We prove in Theorem~\ref
{thm2} that the convergence rate of the fixed-b lag-window estimator is
of order $\log(n)/\sqrt{n}$, better than the fastest rate achievable by
the more commonly used lag-window estimator. Similar comparisons based
on the convergence of $\mathbf{T}_n$ has been reported elsewhere in the
literature. Jansson \cite{jansson04} studied stationary Gaussian moving average
models and established that the rate of convergence of $\mathbf{T}_n$
is $n^{-1}\log(n)$. Sun, Phillips and Jin \cite{sunetal08} obtained the rate $n^{-1}$, under
the main assumption that the underlying process is Gaussian and
stationary. It seems unlikely that the convergence rate $n^{-1}$ will
hold without the Gaussian assumption. However, it is unclear whether
the convergence rate $\log(n)/\sqrt{n}$ obtained in Theorem~\ref{thm2}
is tight.

We organize the paper as follows. Section~\ref{secMC} contains the
main results, including the rate of convergence of the fixed-b
lag-window estimator in Section~\ref{secrate}. We present a simulation
example to illustrate the finite sample properties of the confidence
intervals in Section~\ref{secex}. All the main proofs are postponed to
Section~\ref{secproofs} and the \hyperref[secappendix]{Appendix}.

\subsection{Notation}
Throughout the paper $(\mathsf{X},\mathcal{B})$ denotes a measure
space with a
countably generated sigma-algebra $\mathcal{B}$ with a probability
measure of
interest $\pi$. We denote $L^2(\pi)$ the usual space of
$L^2$-integrable functions with respect to $\pi$,\vspace*{-1.5pt} with norm $\|\cdot\|$
and associated inner product $\anglel \cdot\angler $, and we denote
$L^2_0(\pi)$
the subspace of $L^2(\pi)$ of functions orthogonal to the constants:
$L^2_0(\pi)\stackrel{\mathrm{def}}{=}\{f\in L^2(\pi)\dvt  \int
f(x)\pi(\mathrm{d}x)=0\}$.

For a measurable function $f\dvtx  \mathsf{X}\to\mathbb{R}$, a
probability measure
$\nu$ on $(\mathsf{X},\mathcal{B})$ and a Markov kernel $Q$ on
$\mathsf{X}$, we use the
notation: $\nu(f)\stackrel{\mathrm{def}}{=}\int f(x)\nu(\mathrm
{d}x)$, $\bar f\stackrel{\mathrm{def}}{=}f-\pi(f)$,
$Qf(x)\stackrel{\mathrm{def}}{=}\int f(y)Q(x,\mathrm{d}y)$, and
$Q^jf(x)\stackrel{\mathrm{def}}{=}Q\{Q^{j-1}f\}
(x)$, with $Q^0f(x)=f(x)$. For $V\dvtx  \mathsf{X}\to[0,\infty)$, we define
$\mathcal{L}
_{V}$ as the space of all measurable real-valued functions $f\dvtx \mathsf{X}
\to
\mathbb{R}$ s.t. $|f|_{V}\stackrel{\mathrm{def}}{=}\sup_{x\in
\mathsf{X}}|f(x)|/V(x)<\infty$.
For two\vspace*{-2pt}
probability measures $\nu_1,\nu_2$, we denotes $\|\nu_1-\nu_2\|
_\mathrm{tv}
\stackrel{\mathrm{def}}{=}\sup_{|f|\leq1}|\nu_1(f)-\nu_2(f)|$,
the total variation
distance between $\nu_1$ and $\nu_2$, and $\| \nu_1-\nu_2 \|_{V}
\stackrel{\mathrm{def}}{=}
\sup_{\{f, |f|_V \leq1 \}} |\nu_1(f)-\nu_2(f)|$, its $V$-norm
generalization.

For sequences $\{a_n,b_n\}$ of real nonnegative numbers, the notation
$a_n\lesssim b_n$ means that $a_n\leq c b_n$ for all $n$, and for some
constant $c$ that does not depend on $n$. For a random sequence $\{X_n\}
$, we write $X_n=\mathrm{O}_p(a_n)$ if the sequence $|X_n|/a_n$ is bounded in
probability. We say that $X_n=\mathrm{o}_p(a_n)$ if $X_n/a_n$ converges in
probability to zero as $n\to\infty$.

\section{Monte Carlo confidence intervals for reversible Markov~chains}\label{secMC}
Throughout the paper, $P$ denotes a Markov kernel on $(\mathsf
{X},\mathcal{B})$ that
is reversible with respect to $\pi$. This means that for any pair
$f,g\in L^2(\pi)$, $\anglel f,Pg\angler =\anglel g,Pf\angler $. We
assume that $P$
satisfies the following.

\begin{longlist}[A1]
\item[A1]
For $\pi$-almost all $x\in\mathsf{X}$,
\begin{equation}
\label{ergo} \lim_{n\to\infty}\bigl\|P^n(x,\cdot)-\pi
\bigr\|_\mathrm{tv}=0.
\end{equation}
\end{longlist}

\begin{rem}
Assumption~A1 is very basic. For instance, if $P$ is $\phi$-irreducible, and aperiodic (in addition to being reversible with
respect to $\pi$), then A1 holds. If in addition $P$ is Harris
recurrent, then (\ref{ergo}) holds \textit{for all} $x\in\mathsf
{X}$. If
$P$ is a Metropolis--Hastings kernel, Harris recurrence typically
follows from $\pi$-irreducibility. All these statements can be found,
for instance, in Tierney \cite{tierney94}.
\end{rem}

Throughout the section, unless stated otherwise, $\{X_n, n\geq0\}$ is
a (nonstationary) Markov chain on $(\mathsf{X},\mathcal{B})$ with
transition kernel
$P$ and started at some arbitrary (but fixed) point $x\in\mathsf{X}$ for
which~(\ref{ergo}) holds. The Markov kernel $P$ induces in the usual
way a self-adjoint operator (also denoted $P$) on the Hilbert space
$L_0^2(\pi)$ that maps $h\mapsto Ph$. This\vspace*{-1.5pt} operator $P$ admits a
spectral measure $\mathcal{E}$ on $[-1,1]$, and for $h\in L_0^2(\pi)$
we will
write $\mu_h(\cdot)\stackrel{\mathrm{def}}{=}\anglel h,\mathcal
{E}(\cdot)h\angler $ for the associated
nonnegative Borel measure on $[-1,1]$. Assumption~A1 implies
that $\mu_h$ does not charge $1$ or $-1$, that is $\mu_h(\{-1,1\})=0$.
This is Lemma~5 of Tierney \cite{chanetgeyer94}.

\subsection{Confidence interval for \texorpdfstring{$\pi(h)$}{pi(h)}}
Let $h\in L_0^2(\pi)$. We define
\begin{equation}
\label{sig}
\sigma^2_P(h)\stackrel{\mathrm{def}}{=}
\int_{-1}^1\frac{1+\lambda
}{1-\lambda}\mu _h(
\mathrm{d} \lambda),
\end{equation}
that we call the asymptotic variance of $h$. The terminology comes from
the fact that if the Markov chain is assumed stationary, a calculation
(see, e.g., H{\"a}ggstr{\"o}m and Rosenthal \cite{rosenthalethagg07}, Theorem~4) using the properties of
the spectral measure $\mu_h$ gives
\begin{equation}
\label{limsig}
\lim_{n\to\infty}n\mathbb{E} \Biggl[
\Biggl(n^{-1}\sum_{k=1}^n
h(X_k) \Biggr)^2 \Biggr]=\sigma_P^2(h).
\end{equation}
For nonstationary Markov chains, such as the one considered in this
paper, it is unclear whether~(\ref{limsig}) continues to hold in
complete generality. The estimation of $\sigma_P^2(h)$ is often of
interest because when (\ref{limsig}) holds, $\sigma_P^2(h)/n$
approximates the mean squared error of the Monte Carlo estimate
$n^{-1}\sum_{k=1}^n h(X_k)$. An estimate of $\sigma_P^2(h)$ is often
also sought in order to exploit the Kipnis--Varadhan theorem for
confidence interval purposes. It is known (H{\"a}ggstr{\"o}m and Rosenthal \cite{rosenthalethagg07}, Theorem~4) that $\sigma_P^2(h)$ can also be
written as
\begin{equation}
\label{asympvarP}
\sigma_P^2(h) = \sum
_{\ell=-\infty}^{+\infty}\gamma_{|\ell|}(h),
\end{equation}
where for $\ell\geq0$, $\gamma_\ell(h)\stackrel{\mathrm{def}}{=}\anglel h,P^\ell h\angler $. This
suggests the so-called lag-window estimator of $\sigma_P^2(h)$
\begin{eqnarray}
&& \sigma^2_{b_{n}}\stackrel{\mathrm{def}}{=}\sum
_{\ell
=-n+1}^{n-1}w \biggl(\frac{\ell
}{b_n} \biggr)
\gamma_{n,|\ell|},
\nonumber
\\[-8pt]
\label{lagEst}\\[-8pt]
\eqntext{\displaystyle\mbox{where } \gamma_{n,\ell
}\stackrel{
\mathrm{def}} {=} n^{-1}\sum_{j=1}^{n-\ell}
\bigl(h(X_{j})-\hat\pi_{n}(h) \bigr) \bigl(h(X_{j+\ell})-
\hat\pi_{n}(h) \bigr).}
\end{eqnarray}
In the above display, $\hat\pi_n(h)=n^{-1}\sum_{k=1}^nh(X_k)$,
$1\leq
b_n\leq n$ is an integer such that $b_n\to\infty$, as $n\to\infty$, and
$w \dvtx  \mathbb{R}\to\mathbb{R}$ is an even function ($w(-x)=w(x)$) with support
$[-1,1]$, that is, $w(x)\neq0$ on $(-1,1)$ and $w(x)=0$ for $|x|\geq
1$. Since $w$ has support $[-1,1]$, the actual range for $\ell$ in the
summation defining $\sigma^2_{b_n}$ is $-b_n+1\leq\ell\leq b_n-1$.

The lag-window estimator $\sigma^2_{b_n}$ can be applied more broadly
in time series and the method has a long history. Some of the earlier
work go back to the 1950s (Grenander and Rosenblatt  \cite{grenanderetrosenblatt53}, Parzen \cite{parzen57}).
Convergence results specific to nonstationary Markov chains have been
established recently (see, e.g.,   Damerdji \cite{damerdji95}, Flegal and Jones   \cite{jonesetal09},
Atchad{\'e} \cite{atchadeaos11}  and the references therein);
however, under assumptions that are much stronger than A1. It
remains an open problem whether $\sigma^2_{b_n}$ can be shown to
converge to $\sigma_P^2(h)$ assuming only A1. In particular, the
author is not aware of any result that establishes the consistency of
$\sigma^2_{b_n}$ without assuming that $P$ is geometrically ergodic.

However, if the goal is to construct a confidence interval for $\pi
(h)$, we will now see that it is enough to assume A1 and
$\sigma
_P^2(h)<\infty$. Consider the lag-window estimator obtained by setting
$b_n=n$. This writes
\begin{equation}
\label{lagEst2}
\sigma^2_{n}\stackrel{\mathrm{def}} {=}\sum
_{\ell=-n+1}^{n-1}w \biggl(\frac{\ell
}{n} \biggr)
\gamma_{n,|\ell|}.
\end{equation}
This estimator is well known to be inconsistent for estimating $\sigma
_P^2(h)$, but has recently attracted a lot of interest in the
Econometrics literature under the name of fixed-b asymptotics (Kiefer, Vogelsang and
Bunzel \cite{kieferetal00},
Kiefer and Vogelsang
\cite{kieferetvogelsang05},
Sun,
Phillips and Jin
\cite{sunetal08}, see also Neave \cite{neave70}
for some pioneer work). This paper takes inspiration from this
literature. However, unlike these works, we exploit the Markov
structure and we do not impose any stationary assumption. We introduce
the function $v(t)\stackrel{\mathrm{def}}{=}\int_0^1w(t-u)\,\mathrm{d}u$, $t\in
[0,1]$, and the
kernel $\phi \dvtx   [0,1]\times[0,1]\to\mathbb{R}$, where
\begin{equation}
\label{kernelphi}
\phi(s,t)=w(s-t)-v(s)-v(t) + \int_0^1
v(t)\,\mathrm{d}t,\qquad s,t \in[0,1].
\end{equation}

We say that a kernel $k \dvtx   [0,1]\times[0,1]\to\mathbb{R}$ is positive
definite if for all $n\geq1$, all $a_1,\ldots,a_n\in\mathbb{R}$, and
$t_1,\ldots,t_n\in[0,1]$, $\sum_{i=1}^n\sum_{j=1}^na_ia_j
k(t_i,t_j)\geq0$. We will assume that the weight function $w$ in (\ref
{lagEst2}) is such that the following holds.

\begin{longlist}[A2]
\item[A2]
The function $w \dvtx  \mathbb{R}\to\mathbb{R}$ is an even function, with support
$[-1,1]$, and of class $\mathcal{C}^2$ on $(-1,1)$. Furthermore, the
kernel $\phi$ defined in (\ref{kernelphi}) is positive definite, and
not identically zero.
\end{longlist}

\begin{example}\label{exw2}
Assumption~A2 holds for the function $w$ given by
$w(u)=(1-u^2)\mathbf{1}_{(-1,1)}(u)$. Indeed in this case, a simple
calculation gives that $\phi(s,t)=2(s-0.5)(t-0.5)$, which (by its
multiplicative form) is clearly positive definite. In this particular
case, solving $\int_0^1\phi(s,t)u(t)\,\mathrm{d}t=\alpha u(s)$ yields the unique
eigenvalue $\alpha=2\int_0^1(t-0.5)^2\,\mathrm{d}t=1/6$.
\end{example}

A general approach to guarantee that $\phi$ as in (\ref{kernelphi}) is
positive definite is to start with a positive definite function $w$, as
the next lemma shows.

\begin{lemma}\label{lemposdef}
Suppose that the kernel $[0,1]\times[0,1]\to\mathbb{R}$ defined by
$(s,t)\mapsto w(s-t)$ is continuous and positive definite. Then $\phi$
as in (\ref{kernelphi}) is also positive definite.
\end{lemma}

\begin{pf}
By Mercer's theorem (see Theorem~\ref{theomercer}), there exist
nonnegative numbers $\{\lambda_j, j\geq0\}$, orthonormal functions
$\xi_j \dvtx  [0,1]\to\mathbb{R}$ such that $\int_0^1w(t-s)\xi
_j(s)\,\mathrm{d}s=\lambda
_j\xi_j(t)$, and
\[
w(t-s)=\sum_{j\geq0}\lambda_j
\xi_j(t)\xi_j(s),
\]
and the series converges uniformly and absolutely. It is easy to show
that one can
interchange integral and sum and write $v(t)=\int_0^1w(t-s)\,\mathrm{d}s=\sum_{j\geq0}\lambda_j \xi_j(t)\int_0^1\xi_j(s)\,\mathrm{d}s$,
$\int_0^1v(t)\,\mathrm{d}t=\int_0^1\int_0^1w(t-s)\,
\mathrm{d}s\,\mathrm{d}t=\sum_{j\geq0}\lambda_j (\int_0^1\xi
_j(t)\,\mathrm{d}t )^2$, and then we get
\[
\phi(s,t)=\sum_{j\geq0}\lambda_j \biggl(
\xi_j(t)-\int_0^1\xi
_j(t)\,\mathrm{d}t \biggr) \biggl(\xi_j(s)-\int
_0^1 \xi_j(s)\,\mathrm{d}s \biggr).
\]
This expression of $\phi$ easily shows that it is positive definite.
\end{pf}

The usual approach for showing that the kernel $(s,t)\mapsto w(s-t)$ is
positive definite is by showing that the weight function $t\mapsto
w(t)$ is a characteristic function (or more generally the Fourier
transform of a positive measure) and applying Bochner's theorem. This
approach shows that A2 holds for the Bartlett function
$w(x)=(1-|x|)\mathbf{1}_{(-1,1)}(x)$, the Parzen function
\[
w(x)=\cases{
1-6x^2+6|x|^3,&
 \quad\mbox{if }$|x|\leq\frac{1}{2},$
\vspace*{3pt}\cr
2 (1-|x| )^3, & \quad\mbox{if }$\frac{1}{2}\leq|x|\leq1$,
\vspace*{3pt}\cr
0, &$\quad\mbox{if }|x|>1$,}
\]
and for a number of others weight functions (see, e.g., Hannan \cite{hannan70}, pages 278--279 for details). In the case of the Bartlett
function, the kernel $\phi$ is given by
\[
\label{phibart}
\phi(s,t)=\tfrac{2}{3}-s(1-s)-t(1-t) -|s-t|.
\]
For the Parzen function, we have
\[
v(s)=\frac{3}{8}+s\wedge(1-s) -2 \bigl(s\wedge(1-s) \bigr)^3 +
\bigl(s\wedge(1-s) \bigr)^4\quad \mbox{and}\quad \int_0^1v(t)
\,\mathrm{d}t= \frac{23}{40},
\]
where $a\wedge b\stackrel{\mathrm{def}}{=}\min(a,b)$.
\medskip

Assumption~A2 implies that $\phi$, considered as a linear
operator on $L^2[0,1]$ ($\phi f(s)=\int_0^1 \phi(s,t)f(t)\,\mathrm{d}t$) is
self-adjoint, compact and positive. Therefore, it has only nonnegative
eigenvalues, and a countable number of positive eigenvalues. We denote
$\{\alpha_j,j\in\mathsf{I}\}$ the set of positive eigenvalues of
$\phi$ (each
repeated according to its multiplicity). The index set $\mathsf
{I}\subseteq\{
1,2,\ldots\}$ is either finite or $\mathsf{I}=\{1,2,\ldots\}$. We
introduce the
random variable $\mathbf{T}_w$ defined as
\[
\mathbf{T}_w\stackrel{\mathrm{def}}{=}\frac{Z_0}{\sqrt{\sum_{i\in
\mathsf{I}}\alpha_i
Z_i^2}}\qquad \mbox{where } \{Z_0,Z_i, i\in\mathsf{I}\} \stackrel{\mathrm{i.i.d.}}{\sim}\mathbf{N}(0,1).
\]
Here is the main result.

\begin{theo}\label{thm1}
Assume \textup{A1}--\textup{A2}, and $h\in L^2(\pi)$. If $0<\sigma
_P^2(h)<\infty$, then as $n\to\infty$,
\[
\sigma_n^2\stackrel{\mathsf{w}} {\to}
\sigma_P^2(h)\sum_{i\in
\mathsf
{I}}
\alpha_i Z_i^2\quad \mbox{and}\quad
\mathbf{T}_n\stackrel {\mathrm{def}} {=}\frac
{1}{\sigma_n\sqrt{n}}\sum
_{k=1}^n \bigl(h(X_k)-\pi(h)
\bigr) \stackrel {\mathsf{w}} {\to} \mathbf{T}_w,
\]
where $\{Z_i, i\in\mathsf{I}\}\stackrel{\mathrm{i.i.d.}}{\sim}\mathbf{N}(0,1)$.
\end{theo}

\begin{pf}
See Section~\ref{secproofthm1}.
\end{pf}

The theorem implies that the confidence interval
\begin{equation}
\label{CI2} \hat\pi_n(h)\pm t_{1-\alpha/2}\sqrt{
\frac{\sigma_n^2}{n}},
\end{equation}
is an asymptotically valid Monte Carlo confidence interval for $\pi
(h)$, where $t_{1-\alpha/2}$ is the $(1-\alpha/2)$-quantile of the
distribution of $\mathbf{T}_w$. These quantiles are intractable in
general but can be easily approximated by Monte Carlo simulation (see
Section~\ref{secTw}).

The assumption that $\sigma_P^2(h)$ is finite can be difficult to
check. When $P$ is known to satisfy a drift condition, one can find
whole class of functions for which the asymptotic variance is finite,
as the following proposition shows. The proposition uses Markov chain
concepts that have not been defined above, and we refer the reader to
Meyn and Tweedie \cite{meynettweedie93} for details.

\begin{prop}\label{propMC}
Suppose that $P$ is $\phi$-irreducible and aperiodic, with invariant
distribution $\pi$. Suppose also that there exist measurable functions
$V,f \dvtx  \mathsf{X}\to[1,\infty)$, constant $b<\infty$, and some
petite set
$C\in\mathcal{B}$ such that
\begin{equation}
\label{subgeodrift}
PV(x)\leq V(x)-f(x) +b\mathbf{1}_C(x),\qquad x\in
\mathsf{X}.
\end{equation}
If $\pi(fV)<\infty$, then for all $h\in\mathcal{L}_f$, $\sigma
_P^2(h)<\infty$.
\end{prop}

\begin{pf}
This is a well-known result. We give the proof only for completeness.
Without any loss of generality, suppose that $\pi(h)=0$. We recall that
$\sigma_P^2(h)=\pi(h^2) + 2\sum_{j\geq1}\anglel h, P^j h\angler $.
Since $|\anglel h,P^jh\angler |\leq\int|h(x)||P^jh(x)|\pi(\mathrm{d}x)$, we obtain
\[
\sum_{j\geq0}\bigl| \bigl\anglel h, P^j h \bigr
\angler \bigr|\leq|h|_f\int\bigl|h(x)\bigr| \biggl\{ \sum
_{j\geq
0}\bigl\|P^j(x,\cdot)-\pi(\cdot)\bigr\|_f
\biggr\}\pi(\mathrm{d}x).
\]
Since $P$ is $\phi$-irreducible and aperiodic, and under the drift
condition (\ref{subgeodrift}), Meyn and Tweedie \cite{meynettweedie93}, Theorem~14.0.1
implies that there exists a finite constant $B$ such that $\sum_{j\geq
0}\|P^j(x,\cdot)-\pi(\cdot)\|_f\leq B V(x)$, $x\in\mathsf{X}$. We
conclude that
\[
\sigma_P^2(h)\leq2B|h|_f\int\bigl|h(x)\bigr|V(x)
\pi(\mathrm{d}x)\leq2B|h|^2_f\int f(x)V(x)\pi(
\mathrm{d}x)<\infty.
\]
\upqed\end{pf}

\begin{rem}
Proposition~\ref{propMC} has a number of well-known special cases. The
most common case is when $f=\lambda V$ for some $\lambda\in(0,1)$, in
which case $P$ is geometrically ergodic and $\sigma_P^2(h)<\infty$ for
all $h\in\mathcal{L}_{V^{1/2}}$. Another important special case is
$f=V^{\alpha
}$, for some $\alpha\in[0,1)$. Such drift condition implies that the
Markov chain converges at a polynomial rate. If $\alpha\geq0.5$, then
Proposition~\ref{propMC} implies that $\sigma_P^2(h)<\infty$ for all
$h\in\mathcal{L}_{V^{\alpha-0.5}}$. To see this, notice that (\ref
{subgeodrift}) with $f=V^\alpha$, and Jarner and Roberts \cite{jarneretroberts02}, Lemma~3.5 imply that $PV^{1/2}\leq V^{1/2}-cV^{\alpha-1/2} + b_1\mathbf
{1}_C$. Since $\pi(V^\alpha)<\infty$, the claim follows from
Proposition~\ref{propMC}.
\end{rem}

\subsection{Example: Metropolis Adjusted Langevin Algorithm for smooth~densities}
We give another example where it is possible to check that $\sigma
_P^2(h)<\infty$ without geometric ergodicity. Take $\mathsf
{X}=\mathbb{R}^d$
equipped with the usual Euclidean inner product $\anglel \cdot,\cdot
\angler _2$,
norm $|\cdot|$, and the Lebesgue measure denoted $\mathrm{d}x$. We consider
a probability measure $\pi$ that has a density with respect to the
Lebesgue measure, and in a slight abuse of notation we use the same
symbol to represent $\pi$ and its density: $\pi(x)=\mathrm{e}^{-u(x)}/Z$, for
some function $u \dvtx  \mathsf{X}\to\mathbb{R}$ that we assume is differentiable,
with gradient $\nabla u$.

Let\vspace*{-1pt} $q_\sigma(x,\cdot)$ denotes the density of the Gaussian
distribution $\mathbf{N} (x-\frac{\sigma^2}{2}\rho(x)\nabla
u(x),\sigma^2I_d )$, where the term $\rho(x)\geq0$ is used to
modulate the drift $-\frac{\sigma^2}{2}\nabla u(x)$, and $\sigma>0$ is
a scaling constant. We consider the Metropolis--Hastings algorithm that
generates a Markov chain $\{X_n, n\geq0\}$
with invariant distribution $\pi$ as follows. Given $X_n=x$, we propose
$Y\sim q_\sigma(x,\cdot)$. We either
``accept'' $Y$ and set $X_{n+1}=Y$ with probability $\alpha(x,Y)$, or we
``reject'' $Y$ and set
$X_{n+1}=x$, where
\[
\alpha(x,y)\stackrel{\mathrm{def}} {=}\min \biggl(1,\frac{\pi
(y)}{\pi(x)}
\frac
{q_\sigma
(y,x)}{q_\sigma(x,y)} \biggr).
\]
When $\rho(x)=0$, we get the Random Walk Metropolis (RWM), and when
$\rho(x)=1$, we get the Metropolis Adjusted Langevin Algorithm (MaLa).
However, we are mainly interested in the case where
\begin{equation}
\label{rho}
 \rho(x)\stackrel{\mathrm{def}} {=}\frac{\tau}{\max(\tau,|\nabla
u(x)|)},\qquad x\in
\mathsf{X}
\end{equation}
for some given constant $\tau>0$, which corresponds to the truncated
MaLa proposed by Roberts and Tweedie \cite{robertsettweedie96b}. The truncated MaLa
combines the stability of the RWM and the mixing of the MaLa. It is
known to be geometrically ergodic whenever RWM is geometrically ergodic
(Atchad{\'e} \cite{atchade05}). However, checking in practice that the truncated
MaLa is geometrically ergodic can be difficult, as this involves
checking conditions on the curvature of the log-density. We show in the
next result that if the gradient of the log-density $u$ is Lipschitz
and unbounded then $P$ satisfies a drift condition of the type (\ref
{subgeodrift}), and $\sigma_P^2(h)$ is guaranteed to be finite for
certain functions.

\begin{longlist}[B1]
\item[B1] Suppose that $u$ is bounded from below, continuously
differentiable, and $\nabla u$ is Lipschitz, and
\[
\limsup_{|x|\to\infty}\bigl|\nabla u(x)\bigr|= + \infty.
\]
\end{longlist}

\begin{theo}\label{thmRWM}
Assume \textup{B1} and (\ref{rho}). Set $V(x)\stackrel{\mathrm{def}}{=}a + u(x)$, where
$a\in
\mathbb{R}$ is chosen such that $V\geq1$. Then there exist $b,r\in
(0,\infty
)$ such that
\begin{equation}
PV(x)\leq V(x)-\frac{\sigma^2}{4} \rho(x)\bigl|\nabla u(x)\bigr|^2 +b\mathbf
{1}_{\{|x|\leq r\}}(x),\qquad x\in\mathsf{X}.
\end{equation}
In particular, if $\int u(x)|\nabla u(x)|\mathrm{e}^{-u(x)}\,\mathrm{d}x<\infty$, then
$\sigma_P^2(h)<\infty$ for all $h\in\mathcal{L}_f$, where
$f(x)=\rho
(x)|\nabla u(x)|^2$.
\end{theo}

\begin{pf}
See Section~\ref{secproof-thmRWM}.
\end{pf}

\begin{rem}
This result can be useful in contexts where the $\log$-density $u$ is
known to have a Lipschitz gradient, but is too complicated to allow an
easy verification of the geometric ergodicity conditions.
\end{rem}

\subsection{On the distribution of the random variable $\mathbf{T}_w$}\label{secTw}
It is clear that the limiting distribution $\mathbf{T}_w$ used for
constructing the confidence interval (\ref{CI2}) depends on the choice
of $w$. More research is needed to explain how to best choose $w$ in
this regard. But from the limited simulations done in this paper, we
found that weight functions $w$ with large characteristic exponents
lead to heavy-tailed limiting distributions $\mathbf{T}_w$, and wider
confidence intervals. The characteristic exponent of a weight function
$w$ is the largest number $r>0$ such that $\lim_{u\to
0}|u|^{-r}(1-w(u))\in(0,\infty)$.
Overall, we recommend the use of the
Bartlett weight function $w(u)=(1-|u|)\mathbf{1}_{(-1,1)}(u)$, which
has characteristic exponent $1$, and has behaved very well in the
simulations conducted.

Another issue is how to compute the quantiles of $\mathbf{T}_w$. As
defined, the distribution of $\mathbf{T}_w$ is intractable in general,
as it requires knowing the eigenvalues of $\phi$. But the next result
gives a straightforward method for approximate simulation from $\mathbf{T}_w$.

\begin{prop}\label{proptheo2}
Let $\{Z_j, 1\leq j\leq N\}$ be i.i.d. standard normal random
variables. Then
\[
\mathbf{T}_w^{(N)}\stackrel{\mathrm{def}} {=}
\frac{\sum_{j=1}^N
Z_j}{\sqrt{\sum_{i=1}^N\sum_{j=1}^N\phi (\frac{i-1}{N},\frac{j-1}{N}
)Z_iZ_j}}\stackrel {\mathsf{w}} {\to}\mathbf{T}_w \qquad\mbox{as } N\to\infty.
\]
\end{prop}

\begin{rem}
As pointed out by a referee, one can also approximately sample from
$\mathbf{T}_w$ by generating $X_{1:N}\stackrel{\mathrm{i.i.d.}}{\sim} \mathbf
{N}(0,1)$, and compute $T_N$, with $h(x)=x$. The approach in
Proposition~\ref{proptheo2} is similar, but replaces $\sigma_N^2$ by
$\check\sigma_N^2$ as defined in (\ref{checksigman}). By Lemma~\ref
{lemtechlem1thm1}, the two approaches are essentially equivalent.
\end{rem}

\begin{pf*}{Proof of Proposition~\protect\ref{proptheo2}}
Let $\{0,\alpha_j, j\in\mathsf{I}\}$ be the eigenvalues of $\phi$, with
associated eigenfunctions $\{\Psi_0,\Psi_j, j\in\mathsf{I}\}$
($\Psi
_0\equiv
1$). By Mercer's theorem (see Theorem~\ref{mercer} in the \hyperref[secappendix]{Appendix}),
\[
\sum_{i=1}^N\sum
_{k=1}^N\phi \biggl(\frac{i-1}{N},
\frac
{k-1}{N} \biggr)Z_iZ_k=N\sum
_{j\in\mathsf{I}}\alpha_j \Biggl(\frac{1}{\sqrt
{N}}\sum
_{i=1}^N\Psi _j \biggl(
\frac{i-1}{N} \biggr)Z_i \Biggr)^2.
\]
Hence,
\[
\mathbf{T}_w^{(N)}=\frac{{1}/{\sqrt{N}}\sum_{i=1}^N\Psi
_0 (({i-1})/{N} )Z_i}{\sqrt{\sum_{j\in\mathsf{I}}\alpha_j
({1}/{\sqrt{N}}\sum_{i=1}^N\Psi_j (({i-1})/{N} )Z_i )^2}}.
\]
It is an application of Lemma~\ref{techlem12} that as $N\to\infty$,
$\{
\frac{1}{\sqrt{N}}\sum_{i=1}^N\Psi_0 (\frac{i-1}{N}
)Z_i,\break \frac
{1}{\sqrt{N}}\sum_{i=1}^N\Psi_j (\frac{i-1}{N} )Z_i, 
j\in\mathsf{I}\}
$ converges weakly to $\{Z_0,Z_j, j\in\mathsf{I}\}$. The result then follows
from the continuous mapping theorem.
\end{pf*}

We use Proposition~\ref{proptheo2} to approximately simulate $\mathbf
{T}_w$ for the function $w(u)=(1-u^2)\mathbf{1}_{(-1,1)}(u)$, and for
the Bartlett and Parzen functions. Table~\ref{tab1} reports the $95\%$ and
$97.5\%$ quantiles, computed based on 10\,000 independent samples of
$\mathbf{T}_w^{(N)}$, with $N={}$3000. We replicate these estimates $50$
times to evaluate the Monte Carlo errors reported in parenthesis.

As explained in Example~\ref{exw2}, in the case $w(u)=(1-u^2)\mathbf
{1}_{(-1,1)}(u)$, $\mathbf{T}_w=\sqrt{6}T_1$, where $T_\nu$ denotes the
student's distribution with $\nu$ degree of freedom; thus, is this case
we can compute accurately the quantiles. In particular, the $95\%$ and
$97.5\%$ quantiles are $15.465$ and $31.123$, respectively.

\begin{table}
\tablewidth=210pt
\caption{Approximations of $t$ such that $\mathbb{P}(\mathbf
{T}_w>t)=\alpha/2$}\label{tab1}
\begin{tabular*}{210pt}{@{\extracolsep{\fill}}lll@{}}
\hline
&${\alpha=10\%}$&${\alpha=5\%}$\\
\hline
$w(u)=(1-u^2)_+$ & 15.49 (0.06) & 31.21 (0.19)\\
Parzen & \phantom{0}4.11 (0.01) & \phantom{0}5.64 \\
Bartlett & \phantom{0}3.77 (0.005) & \phantom{0}4.78 (0.01)\\
\hline
\end{tabular*}
\end{table}

\subsection{Rate of convergence of \texorpdfstring{$\sigma_n^2$}{$sigma_n^2$}}\label{secrate}
An interesting question is understanding how the lag-window estimators
$\sigma_n^2$ and $\sigma_{b_n}^2$ compare. On one hand, the asymptotic
behavior of $\sigma_{b_n}^2$ is better understood. In the stationary
case, the best rate of convergence of $\sigma_{b_n}^2$ towards $\sigma
_P^2(h)$ is $n^{-q/(1+2q)}$ (see, e.g., Parzen \cite{parzen57}, Theorem~5A--B),
where $q$ is the largest number $q\in(0,r]$ such that $\sum_{j\geq1}
j^q\gamma_j(h)<\infty$, where $\gamma_j(h)=\anglel h,P^jh\angler $,
and $r$ is
the characteristic exponent of $w$. This optimal rate is achieved by
choosing $b_n\propto n^{{1}/({1+2q})}$. Hence, the optimal rate in
the case of a geometrically ergodic Markov chain is $n^{-r/(1+2r)}$.
However, it is well documented (see, e.g., Newey and West \cite{neweyetwest94}) that
the finite sample properties of $\sigma_{b_n}^2$ are very sensitive to
the actual constant in $b_n\propto n^{{1}/({1+2q})}$, and some tuning
is often required in practice. On the other hand, the fixed-b framework
has the advantage that it requires no tuning, since $b_n=n$.
Furthermore, we establish in this section that $\sigma_n^2$ has a
better convergence rate. Reversibility plays no role in this
discussion. We further simplify the analysis by assuming that $P$
satisfies a geometric ergodicity assumption:
\begin{longlist}[(G)]
\item[(G)] There exists a measurable function $V\dvtx   \mathsf
{X}\to
[1,\infty)$ such that $\pi(V)<\infty$, and for all $\beta\in(0,1]$,
\begin{equation}
\label{rateconv}
\bigl\| P^n(x,\cdot) - \pi(\cdot) \bigr\|_{V^\beta} \leq
C \rho^n \ V^{\beta
}(x),\qquad  n\geq0, x\in\mathsf{X}.
\end{equation}
\end{longlist}

Denote $\mathsf{Lip}_1(\mathbb{R})$ the set of all bounded Lipschitz functions
$f\dvtx  \mathbb{R}\to\mathbb{R}$ such that
\[
|f|_{\mathrm{Lip}}\stackrel{\mathrm{def}} {=}\sup_{x\neq y}
\frac
{|f(x)-f(y)|}{|x-y|}\leq1.
\]
For $P,Q$ two probability measures on $\mathbb{R}$, we define
\[
\mathsf{d}_1(P,Q)\stackrel{\mathrm{def}} {=}\sup
_{f\in\mathsf
{Lip}_1(\mathbb{R})} \biggl\llvert \int f\,\mathrm{d}P-\int f\,\mathrm{d}Q
\biggr\rrvert.
\]
$\mathsf{d}_1(P,Q)$ is the Wasserstein metric between $P,Q$.
An upper bound on $\mathsf{d}_1(P_n,P)$ gives a Berry--Esseen-type
bound on the rate of weak convergence of $P_n$ to $P$.
In a slight abuse of notation, if $X,Y$ are random variables, and
$X\sim P$ and $Y\sim Q$, we shall also write $\mathsf{d}_1(X,Y)$ to
mean $\mathsf{d}_1(P,Q)$.

\begin{theo}\label{thm2}
Suppose that \textup{A2} and \textup{(G)} hold. Suppose
also that $\mathsf{I}$ is finite. For $\delta\in[0,1/4)$, let $h\in
\mathcal{L}
_{V^{\delta}}$ be such that $\pi(h)=0$, and $\sigma_P^2(h)=1$. Then
\begin{equation}
\label{ratetheo2}
\mathsf{d}_1 \bigl(\sigma_n^2,
\chi ^2 \bigr)\lesssim\frac{\log(n)}{\sqrt{n}}\qquad \mbox{as }n\to \infty,
\end{equation}
where $\chi^2=\sum_{i\in\mathsf{I}}\alpha_i Z_i^2$, $\{Z_i, i\in
\mathsf{I}\}$ are
i.i.d. $\mathbf{N}(0,1)$, and $\{\alpha_i, i\in\mathsf{I}\}$ is
the set of
positive eigenvalues of $\phi$.
\end{theo}

\begin{pf}See Section~\ref{proofthm2}.
\end{pf}

\begin{rem}
The assumption that $\mathsf{I}$ is finite is mostly technical and it seems
plausible that this result continues to hold without that assumption.
For example, $\mathsf{I}$ is finite for the kernel
$w(u)=(1-u^2)\mathbf
{1}_{(-1,1)}(u)$.
\end{rem}

\subsection{A simulation example}\label{secex}
This section illustrates the finite sample behavior of the fixed-b
confidence interval procedure. We will compare the fixed-b procedure
and the standard confidence interval procedure based on $\sigma_{b_n}^2$ (using a Gaussian limit). As example, we consider the
posterior distribution of a logistic regression model, and use the
Random Walk Metropolis algorithm (Robert and Casella \cite{robertetcasella04}).

Let $\mathsf{X}=\Theta=\mathbb{R}^d$ equipped with its Borel
sigma-algebra, and
$\pi$ be absolutely continuous w.r.t. the Lebesgue measure $\mathrm{d}
\theta$ with
density still denoted by $\pi$. We write $|\theta|$ for the Euclidean
norm of~$\theta$. Let $q_\Sigma$ denotes the density of the normal
distribution $\mathbf{N}(0,\Sigma)$ on $\Theta$ with covariance matrix
$\Sigma$. The Random Walk Metropolis algorithm (RWMA) is a popular MCMC
algorithm that generates a Markov chain with invariant distribution
$\pi
$ and transition kernel given by
\[
P_\Sigma(\theta,A)=\mathbf{1}_A(\theta) + \int
_\mathsf{X}\alpha (\theta ,\theta+z) \bigl(\mathbf{1}_A(
\theta+z)-\mathbf{1}_A(\theta) \bigr)q_\Sigma(z)\,
\mathrm{d}z,\qquad \theta \in\Theta, A\in\mathcal{B}(\Theta),
\]
where $\mathbf{1}_A$ denotes the indicator function, and $\alpha
(\theta
,\vartheta)\stackrel{\mathrm{def}}{=}\min (1,\frac{\pi
(\vartheta)}{\pi(\theta
)} )$
is the acceptance probability.

We assume that $\pi$ is the posterior distribution from a logistic
regression model. More precisely, we assume that we have binary
responses $y_i\in\{0,1\}$, where
\[
y_i\sim\mathcal{B} \bigl(p \bigl(x_i'
\theta \bigr) \bigr),\qquad  i=1,\ldots,n,
\]
and $x_i\in\mathbb{R}^d$ is a vector of covariate, and $\theta\in
\mathbb{R}^d$
is the vector of parameter. $\mathcal{B}(p)$ denotes the Bernoulli
distribution with parameter $p\in(0,1)$, and $p(x)=\frac{\mathrm{e}^x}{1+\mathrm{e}^x}$
is the cdf of the logistic distribution. Let $X\in\mathbb{R}^{n\times d}$
denote the matrix with $i$th row $x_i'$. Let $\ell(\theta\vert X)$
denotes the log-likelihood function of the model. We assume a Gaussian
prior $\mathbf{N}(0,s^2I_d)$ for $\theta$, with $s=20$. The posterior
distribution of $\theta$ then becomes
\[
\pi (\theta|X )\propto \mathrm{e}^{\ell(\theta\vert
X)}\mathrm{e}^{-{1}/({2s^2})|\theta|^2}.
\]

It is known that for this target distribution the RWM is geometrically
ergodic (see, e.g., Atchad{\'e} \cite{atchadeaos11}, Section~5.2). Therefore, for
all polynomial functions Theorem~\ref{thm1} holds. It is also known
that with an appropriate choice of $b_n$, $\sigma^2_{b_n}$ converges in
probability to $\sigma_P^2(h)$ (see, e.g., Atchad{\'e} \cite{atchadeaos11}, Theorem~4.1, and Corollary~4.1). So we will compare the fixed-b confidence
intervals and the classical confidence intervals based on $\sigma_{b_n}^2$.

We simulate a Gaussian dataset with $n=250$, $d=15$, and simulate the
components of the true value of $\beta$ from a $\mathsf{U}(-10,10)$. We
first run the adaptive chain for $10^6$ iterations and take the sample
posterior mean of $\beta$ as the ``true'' posterior mean. We focus on
the coefficient $\beta_1$. Each sampler is run for 30\,000 iterations,
with no burn-in period. For the RMW, we use a covariance matrix $\Sigma
=c \mathrm{I}_{15}$, where $c$ is chosen such that the acceptance
probability in stationarity is about $30\%$, obtained from a
preliminary run.

\begin{figure}

\includegraphics{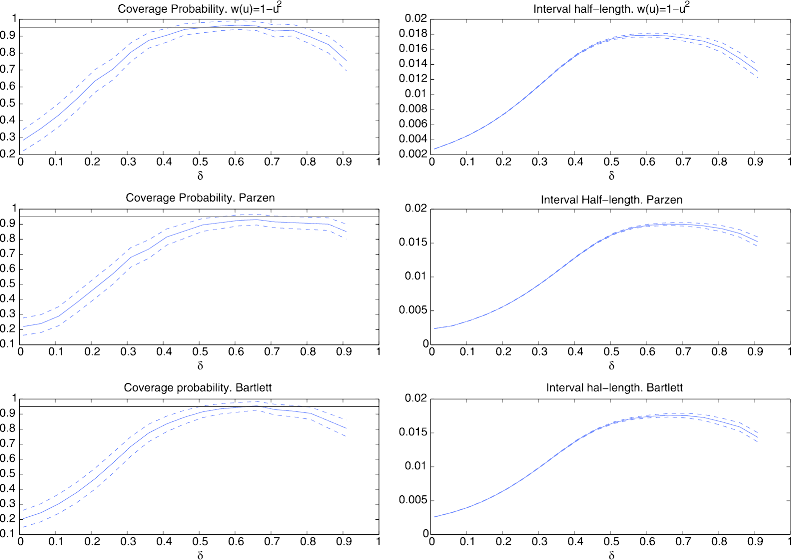}

\caption{Coverage probability and confidence
interval half-length for parameter $\beta_1$ for different values of~$\delta$ using $\sigma_{b_n}^2$, and $b_n=n^{\delta}$. The dashed line
is the $95\%$ confidence band estimated from $200$ replications.}
\label{fig1}
\end{figure}

From each sampler, we compute the fixed-b $95\%$ confidence interval,
and a classical $95\%$ confidence interval. To explore the range of
behavior of the classical procedure, we use $b_n=n^{\delta}$ for
different values of $\delta\in(0,1)$. To estimate coverage probability
and half-length of these confidence intervals, $K=200$ replications are
performed. The result is summarized in Table~\ref{tab2} for the fixed-b
procedure, and in Figure~\ref{fig1} for the classical procedure.

We see from the results that using $b_n=n$ gives very good coverage,
except for the choice $w(u)=(1-u^2)_+$, which generates significantly
wider intervals. This is somewhat expected given the very heavy tail of
the limiting distribution. The result also shows that the confidence
interval procedure based on $\sigma^2_{b_n}$ works equally well when
$b_n$ is carefully chosen, but can perform poorly otherwise.

\begin{table}
\tablewidth=200pt
\caption{Coverage probability and half-length for fixed-b confidence
intervals}\label{tab2}
\begin{tabular*}{\tablewidth}{@{\extracolsep{\fill}}lll@{}}
\hline
&{Coverage}& {Half-length}\\
\hline
$w(u)=(1-u^2)_+$ & $0.945 \pm 0.03$ & $0.10 \pm 0.01$\\
Parzen & $0.94 \pm  0.03$ & $0.03 \pm 0.002$\\
Bartlett & $0.955 \pm 0.03$ & $0.02 \pm 0.001$\\
\hline
\end{tabular*} \vspace*{-12pt}
\end{table}

\begin{figure}

\includegraphics{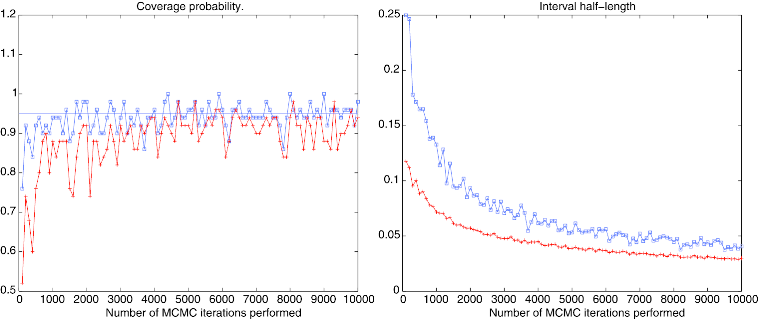}

\caption{Coverage probability and confidence
interval half-length for parameter $\beta_1$ as function of number of
MCMC iterations. The square-line corresponds to using $\sigma^2_n$.}
\label{fig2}
\end{figure}

We also test the conclusion of Theorem~\ref{thm2} by comparing the
finite sample convergence rate of the two confidence interval
procedures. Here, we use only the Bartlett function. For the standard
procedure, we use the best choice of $\delta$ ($\delta\approx0.66$),
as given by the previous simulation. We compute the confidence
intervals after MCMC runs of length $n$, where $n\in\{100,\ldots
,10^4\}$. Each run is repeated $30$ times to approximate the coverage
probabilities and interval lengths. The result is plotted on Figure~\ref{fig2},
and is consistent with Theorem~\ref{thm2} that the fixed-b procedure
has faster convergence. The price to pay is a (slightly) wider interval
length as seen on Figure~\ref{fig2}.

\section{Proofs}\label{secproofs}
\subsection{Proof of Theorem~\texorpdfstring{\protect\ref{thm1}}{2.2}}\label{secproofthm1}
Let $\phi$ as in (\ref{kernelphi}). Assumption~A2 and Mercer's
theorem implies that the kernel $\phi$ has a countable number of
positive eigenvalues $\{\alpha_i, i\in\mathsf{I}\}$ with associated
eigenfunctions $\{\Psi_j, j\in\mathsf{I}\}$ such that
\begin{equation}
\label{mercer}
\phi(s,t)=\sum_{j\in\mathsf{I}}
\alpha_j\Psi_j(s)\Psi_j(t),\qquad (s,t)\in [0,1]
\times[0,1],
\end{equation}
where the convergence of the series is uniform on $[0,1]\times[0,1]$.
Since $\int_0^1\phi(s,t)\,\mathrm{d}t=0$, $0$ is also an eigenvalue of $\phi$ with
eigenfunction $\Psi_0(x)\equiv1$. Hence, we define $\bar\mathsf
{I}=\{0\}
\cup\mathsf{I}
$, $\alpha=\{\alpha_j, j\in\bar\mathsf{I}\}$, with $\alpha_0=0$, and
$\ell
^2(\alpha)$ the associated Hilbert space of real numbers sequences $\{
x_j, j\in\bar\mathsf{I}\}$ such that $\sum_jx_j^2<\infty$,
equipped with the
norm $\|x\|_\alpha=\sqrt{\sum_j\alpha_jx_j^2}$ and the inner product
$\anglel x,y\angler _\alpha\stackrel{\mathrm{def}}{=}\sum_j\alpha
_jx_jy_j$. We will need the
differentiability of the eigenfunction $\Psi_j$. This is given by
Kadota's theorem (Kadota \cite{kadota67}). Under the assumption that $w$ is
continuously twice differentiable, the eigenfunctions $\Psi_j$, $j\in
\mathsf{I}
$ are continuously differentiable (with derivative $\Psi'$) and
\begin{equation}
\label{kadota}
\frac{\partial^2}{\partial s\partial t}\phi(s,t)=\sum_{j\in\mathsf{I}
}
\alpha _j\Psi'_j(s)\Psi'_j(t),
\qquad (s,t)\in[0,1]\times[0,1],
\end{equation}
where again the convergence of the series is uniform on $[0,1]\times
[0,1]$. The expansions (\ref{mercer}) and~(\ref{kadota}) easily imply that
\begin{eqnarray}
 &&\sum_{j\in\mathsf{I}} \alpha_j <
\infty, \qquad\sup_{t\in
[0,1]}\sum_{j\in
\mathsf{I}}
\alpha_j\bigl|\Psi_j(t)\bigr|^2<\infty \quad \mbox{and}
\nonumber
\\[-8pt]
\label{boundPsi}
\\[-8pt]
\nonumber
&& \sup_{t\in
[0,1]}\sum_{j\in\mathsf{I}}
\alpha_j\bigl|\Psi_j'(t)\bigr|^2 < \infty.
\end{eqnarray}

It is easy to check that $\sigma_n^2$ can also be written as
\begin{eqnarray*}
\sigma^2_n & = & \frac{1}{n}\sum
_{i=1}^n \sum_{j=1}^nw
\biggl(\frac{i-j}{n} \biggr) \bigl(\bar h(X_i)-
\pi_n(\bar h) \bigr) \bigl(\bar h(X_{j})-\pi_n(
\bar h) \bigr)
\\
& =& \frac{1}{n}\sum_{i=1}^n\sum
_{j=1}^n \biggl\{w \biggl(
\frac
{i-j}{n} \biggr) -v_{n,i}-v_{n,j}+u_n
\biggr\}\bar h(X_i)\bar h(X_j),
\end{eqnarray*}
where $v_{n,i}=n^{-1}\sum_{\ell=1}^n w (\frac{i-\ell}{n} )$,
and $u_n=n^{-2}\sum_{i=1}^n \sum_{=1}^n w (\frac{i-j}{n} )$.
Notice that $v_{n,i}$ is a Riemann sum approximation of $v(i/n)$, where
$v(t)\stackrel{\mathrm{def}}{=}\int_0^1w(t-u)\,\mathrm{d}u$, and $u_n$
approximates $\int_0^1\!\int_0^1
w(t-u)\,\mathrm{d}u\,\mathrm{d}t=\int_0^1v(t)\,\mathrm{d}t$. In view of this, we introduce
\begin{eqnarray}
\check\sigma_n^2 & \stackrel{
\mathrm{def}} {=}& \frac{1}{n}\sum_{i=1}^n
\sum_{j=1}^n \biggl\{ w \biggl(
\frac{i-j}{n} \biggr) -v \biggl(\frac{i-1}{n} \biggr)-v \biggl(
\frac
{j-1}{n} \biggr)+\int_0^1v(t)\,
\mathrm{d}t \biggr\}\bar h(X_i)\bar h(X_j)
\nonumber
\\[-8pt]
\label{checksigman}
\\[-8pt]
\nonumber
&=& \frac{1}{n}\sum_{i=1}^n\sum
_{j=1}^n\phi \biggl(\frac
{i-1}{n},
\frac
{j-1}{n} \biggr)\bar h(X_i)\bar h(X_j)=\sum
_{\ell\in\mathsf
{I}}\alpha_\ell \Biggl(\frac{1}{\sqrt{n}}
\sum_{i=1}^n\Psi_\ell \biggl(
\frac
{i-1}{n} \biggr)\bar h(X_i) \Biggr)^2.\quad
\end{eqnarray}
The last equality uses the Mercer's expansion for $\phi$ as given in
(\ref{mercer}). This implies that
\begin{eqnarray*}
\mathbf{T}_n &=& \frac{\sum_{i=1}^n \bar h(X_i)}{\sigma_n\sqrt
{n}}\\
&=& \frac{{1}/({\sigma_P(h)\sqrt{n})}\sum_{i=1}^n \Psi_0
(({i-1})/{n} )\bar h(X_i)}{\sqrt{\sum_{\ell\in\mathsf{I}}\alpha
_\ell
({1}/({\sigma_P(h)\sqrt{n})}\sum_{i=1}^n\Psi_\ell (({i-1})/{n} )\bar h(X_i) )^2 +
{(\sigma_n^2-\check\sigma_n^2)}/{\sigma^2_P(h)}}}.
\end{eqnarray*}

Hence, the proof of the theorem boils down to the limiting behavior of
the $\ell^2(\alpha)$-valued process
\[
\Biggl\{ \frac{1}{\sigma_P(h)\sqrt{n}}\sum_{i=1}^n
\Psi_j \biggl(\frac
{i-1}{n} \biggr)\bar h(X_i), j\in
\bar\mathsf{I} \Biggr\},
\]
and the remainder $(\sigma_n^2-\check\sigma_n^2)$. In Lemma~\ref
{lemtechlem2thm1}, we show that $ \{ \frac{1}{\sigma_P(h)\sqrt
{n}}\sum_{i=1}^n\Psi_\ell (\frac{i-1}{n} )\bar h(X_i),\break
\ell
\in\bar\mathsf{I} \}$ converges weakly to $\{Z_\ell, \ell
\in\bar
\mathsf{I}\}$,
and that $\sigma_n^2-\check\sigma_n^2$ converges in probability to
zero. This is done first in the stationary case in Lemmas \ref{techlem12}--\ref{lemtechlem1thm1}, and in the nonstationary case in
Lemma~\ref{lemtechlem2thm1}. Hence, the theorem follows by applying
Slutszy's theorem and the continuous mapping theorem. Everything rely
on a refinement of the martingale approximation of Kipnis and Varadhan \cite{kv86} that we
establish first in Lemma~\ref{techlem1}.

\subsubsection{Martingale approximation for Markov chains}\label{MApproxMC}
Throughout this section, unless stated otherwise, $\{X_n, n\geq0\}$
denotes a stationary reversible Markov chain with invariant
distribution $\pi$ and transition kernel $P$, and we fix $h\in
L_0^2(\pi
)$. We denote $\mathcal{F}_n\stackrel{\mathrm{def}}{=}\sigma
(X_0,\ldots,X_n)$. We introduce the
probability measure $\bar\pi(\mathrm{d}x,\mathrm{d}y)=\pi(\mathrm{d}x)P(x,\mathrm{d}y)$ on $\mathsf
{X}\times
\mathsf{X}$, and we denote $L^2(\bar\pi)$ the associated $L^2$-space with
norm $\nnorm{f}^2\stackrel{\mathrm{def}}{=}\int\!\!\int|f(x,y)|^2\times\pi
(\mathrm{d}x)P(x,\mathrm{d}y)$. For
$\varepsilon>0$, define
\[
U_\varepsilon(x)\stackrel{\mathrm{def}} {=}\sum
_{j\geq0}\frac
{1}{(1+\varepsilon)^{j+1}} P^j h(x),\qquad
G_\varepsilon(x,y)\stackrel{\mathrm{def}} {=}U_\varepsilon
(y)-PU_\varepsilon(x).
\]
Since $P$ is a contraction of $L_0^2(\pi)$, it is clear that
$U_\varepsilon
\in L^2(\pi)$, and $G_\varepsilon\in L^2(\bar\pi)$. Furthermore, for all
$\varepsilon>0$,
\begin{equation}
\label{normU}
\|U_\varepsilon\|\leq\varepsilon^{-1}\|h\| \quad \mbox{and}\quad \nnorm {G_\varepsilon}\leq2\|U_\varepsilon\|.
\end{equation}

When $\sigma_P^2(h)<\infty$ a stronger conclusion is possible, and this
is the key observation made by Kipnis and Varadhan \cite{kv86}, Theorem~1.3. We summarize
their result as follows.

\begin{lemma}[(Kipnis and Varadhan \protect\cite{kv86})] \label{theokv}
Suppose that $h\in L_0^2(\pi)$, and $\sigma_P^2(h)<\infty$. Then for
any sequence $\{\varepsilon_n, n\geq0\}$ of positive numbers such that
$\lim_n\varepsilon_n=0$,
\[
\lim_{n\to\infty} \sqrt{\varepsilon_n}\|U_{\varepsilon_n}
\| =0.
\]
Furthermore, there exists $G\in L^2(\bar\pi)$, with $\int
P(x,\mathrm{d}z)G(x,z)=0$ ($\pi$-a.e.) such that $\sigma_P^2(h)=\nnorm{G}^2$,
and $\lim_{n} \nnorm{G_{\varepsilon_n}-G}=0$.
\end{lemma}

For $n\geq1$, define the process
\[
B_n(t)=\frac{1}{\sigma_P(h)\sqrt{n}}\sum_{i=1}^{\lfloor nt\rfloor
}G(X_i,X_{i-1}),
\qquad 0\leq t\leq1,
\]
and let $\{B(t), 0\leq t\leq1\}$ denotes the standard Brownian
motion. It is an easy consequence of Lemma~\ref{theokv} that $\{
G(X_i,X_{i-1}), 1\leq i\leq n\}$ is a stationary martingale difference
sequence with finite variance. Therefore, by the weak invariance
principle for stationary martingales, $B_n\stackrel{\mathsf{w}}{\to} B$
in $D[0,1]$ equipped with the Skorohod metric. In Corollary~1.5, \cite{kv86}, it is shown that the Markov chain $\{X_n, n\geq0\}$ inherits
this weak invariance principle. For the purpose of this paper, we need
some refinements of this result. Let $\{a_{n,k}, 0\leq k\leq n\}$ be a
sequence of real numbers. Set $|a_n|_\infty\stackrel{\mathrm{def}}{=}\sup_{0\leq k\leq
n}|a_{n,k}|$, and $|a_n|_\mathrm{tv}\stackrel{\mathrm{def}}{=}\sum_{k=1}^n |a_{n,k}-a_{n,k-1}|$.

\begin{lemma}\label{techlem1}
Let $h\in L_0^2(\pi)$ be such that $\sigma_P^2(h)<\infty$.
\begin{enumerate}[(2)]
\item[(1)] If $|a_n|_\infty+|a_n|_\mathrm{tv}$ is bounded in $n$, then
\begin{equation}
\sum_{i=1}^n a_{n,i-1}h(X_i)=
\sum_{i=1}^n a_{n,i-1}G(X_i,X_{i-1})
+ R_n,
\end{equation}
where $n^{-1}\mathbb{E}(|R_n|^2)\to0$ as $ n\to\infty$.
\item[(2)] If $f \dvtx  [0,1]\to\mathbb{R}$ is a continuously differentiable function,
then $\frac{1}{\sigma_P(h)\sqrt{n}}\sum_{i=1}^n f (\frac
{i-1}{n} )\*h(X_i)$ converges weakly to $\int_0^1 f(t)\,\mathrm{d}B(t)$, as
$n\to\infty$.
\end{enumerate}
\end{lemma}

\begin{pf}
Set $S_n\stackrel{\mathrm{def}}{=}\sum_{i=1}^n
a_{n,i-1}h(X_i)$. The function
$U_\varepsilon$ satisfies $(1+\varepsilon) U_\varepsilon(x)-PU_\varepsilon
(x)=h(x)$, $\pi$-a.e. $x\in\mathsf{X}$. This is used to write
\begin{eqnarray*}
a_{n,k-1}h(X_k) &=& a_{n,k-1} \bigl(\varepsilon
U_\varepsilon(X_k) +U_\varepsilon (X_k) -
PU_\varepsilon(X_k) \bigr)\\
&=& a_{n,k-1}\varepsilon
U_\varepsilon(X_k) +a_{n,k-1} \bigl(U_\varepsilon(X_k)-PU_\varepsilon(X_{k-1})
\bigr)
\\
&&{}+ \bigl(a_{n,k-1}PU_\varepsilon(X_{k-1})-a_{n,k}PU_\varepsilon
(X_{k}) \bigr) + (a_{n,k}-a_{n,k-1} )
PU_\varepsilon(X_k).
\end{eqnarray*}
It follows that
\begin{eqnarray*}
S_n &=& \varepsilon\sum_{k=1}^n
a_{n,k-1} U_\varepsilon(X_k) + \sum
_{k=1}^n a_{n,k-1} G(X_k,X_{k-1})
\\
&&{}+ \sum_{k=1}^n a_{n,k-1}
\bigl(G_\varepsilon (X_k,X_{k-1}) -G(X_k,X_{k-1})
\bigr)
\\
&&{}+ \bigl(a_{n,0}PU_\varepsilon(X_0)-a_{n,n}PU_\varepsilon(X_n)
\bigr)+\sum_{k=1}^n (a_{n,k}-a_{n,k-1})PU_\varepsilon(X_{k}),
\end{eqnarray*}
which is valid for any $\varepsilon>0$. In particular with $\varepsilon
=\varepsilon_n=1/n$, we have
\[
S_n=\sum_{k=1}^n
a_{n,k-1} G(X_k,X_{k-1}) + \sum
_{k=1}^n a_{n,k-1} \bigl(G_{\varepsilon_n}(X_k,X_{k-1})
-G(X_k,X_{k-1}) \bigr)+ R_n^{(1)}
+R_n^{(2)} +R_n^{(3)},
\]
where
\begin{eqnarray*}
R_n^{(1)} &\stackrel{\mathrm{def}}{=} & \varepsilon_n
\sum_{k=1}^n a_{n,k-1}
U_{\varepsilon
_n}(X_k),\qquad  R_n^{(2)}\stackrel{\mathrm{def}} {=} \bigl(a_{n,0}PU_{\varepsilon_n} (X_0)-a_{n,n}PU_{\varepsilon_n}(X_n)
\bigr)\quad\mbox{and }
\\
R_n^{(3)} &\stackrel{\mathrm{def}}{=} & \sum
_{k=1}^n (a_{n,k}-a_{n,k-1})PU_{\varepsilon_n}(X_{k}).
\end{eqnarray*}
By stationarity and the martingale property,
\[
\frac{1}{n}\mathbb{E} \Biggl[ \Biggl(\sum_{k=1}^n
a_{n,k-1} \bigl(G_\varepsilon (X_k,X_{k-1})
-G(X_k,X_{k-1}) \bigr) \Biggr)^2 \Biggr]=
\nnorm {G_{\varepsilon
_n}-G}^2\frac{1}{n}\sum
_{k=1}^n a_{n,k-1}^2\to0,
\]
using Lemma~\ref{theokv}, and the assumption on $a_n$. The other
remainders are also easily dealt with.
\[
\frac{1}{\sqrt{n}}\mathbb{E}^{1/2} \bigl(\bigl|R_n^{(3)}\bigr|^2
\bigr)\leq \frac
{1}{\sqrt
{n}}\sum_{k=1}^n
|a_{n,k}-a_{n,k-1}|\mathbb{E}^{1/2}
\bigl(\bigl|PU_{\varepsilon
_n}(X_{k})\bigr|^2 \bigr)=\sqrt{
\varepsilon_n}\|U_{\varepsilon_n}\| |a_n|_\mathrm{tv}
\to0,
\]
using Lemma~\ref{theokv} and the assumption on $a_n$. Similarly,
\begin{eqnarray*}
\frac{1}{\sqrt{n}}\mathbb{E}^{1/2} \bigl(\bigl|R_n^{(2)}\bigr|^2
\bigr) &\leq & 2|a_n|_\infty \sqrt{\varepsilon_n}
\|U_{\varepsilon_n}\|  \to  0\quad \mbox{and}\\
 \frac{1}{\sqrt{n}}\mathbb{E}^{1/2}
\bigl(\bigl|R_n^{(1)}\bigr|^2 \bigr)&\leq & \sqrt {
\varepsilon _n}\|U_{\varepsilon_n}\|\frac{1}{n}\sum
_{k=1}^n |a_{n,k-1}|  \to  0.
\end{eqnarray*}
This proves part (1) of the lemma. For part (2), we use part (1) with
$a_{n,i}=f(i/n)$ to conclude that
\begin{eqnarray*}
&& \frac{1}{\sigma_P(h)\sqrt{n}}\sum_{i=1}^n f \biggl(
\frac
{i-1}{n} \biggr)h(X_i)\\
 &&\quad =  \frac{1}{\sigma_P(h)\sqrt{n}} \sum
_{i=1}^n f \biggl(\frac
{i-1}{n}
\biggr) G(X_i,X_{i-1}) + \mathrm{o}_p(1)
\\
&&\quad=\int_0^1f(t)\,\mathrm{d}B_n(t)
+ \mathrm{o}_p(1),
\end{eqnarray*}
where $A_n=\mathrm{o}_p(1)$ means that $A_n$ converges in probability to zero as
$n\to\infty$. To conclude the proof, it suffices to show that $\int_0^1f(t)\,\mathrm{d}B_n(t)$ converges weakly to $\int_0^1f(t)\,\mathrm{d}B(t)$. This
follows from the weak convergence continuous mapping theorem by
noticing that $B$ has continuous sample path (almost surely), and the
map $D[0,1]\to\mathbb{R}$, $x\mapsto\int_0^1 f(t)\,\mathrm{d}x(t)$
is continuous
at all points $x_0\in C[0,1]$, where the integral $\int_0^1
f(t)\,\mathrm{d}x(t)$ is understood as a Riemann--Stietjes integral. To see the
continuity, take $\{x_n\}$ a sequence of elements in $D[0,1]$ that
converges to $x_0\in C[0,1]$ in the Skorohod metric. Since $x_0\in
C[0,1]$, the sequence $\{x_n\}$ converges to $x_0$ in $C[0,1]$ as well.
By integration by part, $\int_0^1f(t)\,\mathrm{d}x_n(t)=f(1)x_n(1)-f(0)x_n(0)
- \int_0^1 x_n(t)f'(t)\,\mathrm{d}t$, and
\[
\biggl\llvert \int_0^1f(t) \,
\mathrm{d}x_n(t)-\int_0^1f(t)\,\mathrm
{d}x_0(t) \biggr\rrvert \leq |x_n-x_0|_\infty
\biggl(2|f|_\infty+\int_0^1\bigl|f'(t)\bigr|
\,\mathrm{d}t \biggr)\to0,
\]
as $n\to\infty$.
\end{pf}

\begin{lemma}\label{techlem12}
Let $h\in L_0^2(\pi)$ be such that $\sigma_P^2(h)<\infty$. Define
\[
Z^{(n)}\stackrel{\mathrm{def}} {=} \Biggl\{ \frac{1}{\sigma_P(h)\sqrt
{n}}\sum
_{i=1}^n\Psi _j \biggl(
\frac{i-1}{n} \biggr)\bar h(X_i), j\in\bar\mathsf {I} \Biggr\}
\quad \mbox{and}\quad Z\stackrel{\mathrm{def}} {=} \biggl\{\int_0^1
\Psi _j(t)\,\mathrm{d}B(t), j\in\bar\mathsf{I} \biggr\}.
\]
Then as $n\to\infty$, $Z^{(n)}$ converges weakly to $Z$ in $\ell
^2(\alpha)$.
\end{lemma}

\begin{pf}
We need to show that for all $u\in\ell^2(\alpha)$, $\anglel Z^{(n)},u\angler _\alpha\stackrel{\mathsf{w}}{\to}\anglel Z,u\angler _\alpha$,
and that
$\{Z^{(n)}\}$ is tight.

For $u\in\ell^2(\alpha)$, $\anglel Z^{(n)},u\angler _\alpha=\frac
{1}{\sigma
_P(h)\sqrt{n}}\sum_{i=1}^n f_u (\frac{i-1}{n} )\bar h(X_i)$,
where $f_u(t)=\sum_{j}\alpha_j u_j\Psi_j(t)$. From basic results in
calculus, it follows from Kadota's theorem that $f_u$ is continuously
differentiable on $[0,1]$. Hence, by Lemma~\ref{techlem1}, part~(2),
$\anglel Z^{(n)},u\angler _\alpha\stackrel{\mathsf{w}}{\to} \int_0^1f_u(t)\,\mathrm{d}B(t)=\anglel u,Z\angler _\alpha$. To show that $\{Z^{(n)}\}$
is tight,
it suffices to show that
\begin{equation}
\label{tightness} \lim_{N\to\infty}\sup_{n\geq1}
\mathbb{E} \Biggl(\sum_{j=N}^\infty \bigl
\anglel Z^{(n)},e_j \bigr\angler ^2_\alpha
\Biggr)=0.
\end{equation}
We have
\begin{eqnarray*}
\mathbb{E} \bigl( \bigl\anglel Z^{(n)},e_j \bigr\angler
^2_\alpha \bigr)&=& \frac
{\alpha
_j}{\sigma
^2_P(h)n}\sum
_{i=1}^n\sum_{k=1}^n
\Psi_j \biggl(\frac{i-1}{n} \biggr)\Psi _j \biggl(
\frac{k-1}{n} \biggr) \pi \bigl(hP^{|i-k|}h \bigr)
\\
&=& \frac{\alpha_j}{\sigma^2_P(h)n}\int_{-1}^1\sum
_{i=1}^n\sum_{k=1}^n
\Psi _j \biggl(\frac{i-1}{n} \biggr)\Psi_j \biggl(
\frac{k-1}{n} \biggr)\lambda ^{|i-k|}\mu_{h}(\mathrm{d}
\lambda).
\end{eqnarray*}
By Fubini's theorem, for $N\geq1$,
\[
\mathbb{E} \Biggl(\sum_{j=N}^\infty \bigl
\anglel Z^{(n)},e_j \bigr\angler ^2_\alpha
\Biggr)=\frac
{1}{\sigma^2_P(h)n}\int_{-1}^1\sum
_{i=1}^n\sum_{k=1}^n
\sum_{j=N}^\infty \alpha_j
\Psi_j \biggl(\frac{i-1}{n} \biggr)\Psi_j \biggl(
\frac
{k-1}{n} \biggr)\lambda^{|i-k|}\mu_{h}(\mathrm{d}
\lambda).
\]
Let $\varepsilon>0$. By uniform\vspace*{1pt} convergence of the series $\sum_{j}\alpha
_j\Psi_j(s)\Psi_j(t)$, we can find $N_0$ such that for any $N\geq N_0$
and for all $s,t\in[0,1]$, $|\sum_{\ell\geq N}\alpha_\ell\Psi
_\ell
(t)\Psi_\ell(s)|\leq\varepsilon$. So that for all $n\geq1$,
\[
\mathbb{E} \Biggl(\sum_{\ell=N}^\infty \bigl
\anglel Z^{(n)},e_\ell \bigr\angler ^2_\alpha
\Biggr)\leq \frac{\varepsilon}{\sigma^2_P(h)n}\int_{-1}^1\sum
_{i=1}^n\sum_{j=1}^n
\lambda^{|i-j|}\mu_{h}(\mathrm{d}\lambda)\leq
\frac{\varepsilon
}{\sigma
^2_P(h)}\int_{-1}^1\frac{1+\lambda}{1-\lambda}
\mu_{h}(\mathrm{d} \lambda )=\varepsilon,
\]
since $\varepsilon>0$ is arbitrary, this proves (\ref{tightness}).
\end{pf}

\begin{lemma}\label{lemtechlem1thm1}
Let $h\in L_0^2(\pi)$ be such that $\sigma_P^2(h)<\infty$. Then as
$n\to
\infty$, $\mathbb{E} (|\sigma_n^2-\check\sigma_n^2| ) = \mathrm{O}(1/n)$.
Hence $\sigma_n^2-\check\sigma_n^2$ converges in probability to $0$, as
$n\to\infty$.
\end{lemma}

\begin{pf}
Comparing the expression of $\sigma_n^2$ and $\check\sigma_n^2$, we
see that
\begin{eqnarray}
 \sigma_n^2-\check\sigma_n^2 &=&
\biggl(u_n-\int_0^1v(t)\,
\mathrm{d}t \biggr) \Biggl(\frac
{1}{\sqrt{n}}\sum_{i=1}^n
\bar h(X_i) \Biggr)^2
\nonumber
\\[-8pt]
\label{eqdiffsigma}
\\[-8pt]
\nonumber
&&{}- 2 \Biggl(\frac{1}{\sqrt{n}}\sum_{i=1}^n
\bar h(X_i) \Biggr) \Biggl(\frac
{1}{\sqrt{n}}\sum
_{i=1}^n \biggl( v_{n,i}-v \biggl(
\frac{i-1}{n} \biggr) \biggr)\bar h(X_i) \Biggr).
\end{eqnarray}
Since the sequence $\mathbb{E} [ (\frac{1}{\sqrt{n}}\sum_{i=1}^n
\bar
h(X_i) )^2 ]$ converges to the finite limit $\sigma^2(h)$ by
assumption, it is bounded, and there exists a finite constant $c_1$
such that
\begin{eqnarray*}
&& \mathbb{E} \bigl(\bigl|\sigma_n^2-\check\sigma_n^2\bigr|
\bigr)\\
&&\quad\leq   c_1^2 \biggl\llvert u_n-\int
_0^1v(t)\,\mathrm{d}t \biggr\rrvert +
\frac{2c_1}{n} \mathbb{E}^{1/2} \Biggl[ \Biggl(\frac
{1}{\sqrt
{n}}\sum
_{i=1}^n n \biggl(v_{n,i}-v
\biggl( \frac{i-1}{n} \biggr) \biggr) h(X_i)
\Biggr)^2 \Biggr].
\end{eqnarray*}
Set $a_{n,0}=0$, $a_{n,i}\stackrel{\mathrm{def}}{=}n
(v_{n,i}-v (\frac
{i-1}{n} ) )$. We recall that $v_{n,i}=n^{-1}\sum_{\ell=1}^n
w (\frac{i-\ell}{n} )$, and $v(t)=\int_0^1w(t-u)\,\mathrm{d}u$,
and write
\begin{eqnarray*}
a_{n,i}&=&n\sum_{\ell=1}^n \int
_{(\ell-1)/n}^{\ell/n} \biggl[w \biggl(\frac
{i-1}{n}-
\frac{\ell-1}{n} \biggr) - w \biggl(\frac{i-1}{n}-u \biggr) \biggr] \,
\mathrm{d}u
\\
&=& n\sum_{\ell=1}^n \int
_{(\ell-1)/n}^{\ell/n} \biggl(\frac{\ell
-1}{n}-u \biggr)\int
_0^1 w' \biggl(\frac{i-1}{n}-
\frac{\ell
-1}{n}-t \biggl(u-\frac{\ell-1}{n} \biggr) \biggr)\,\mathrm{d}t \,
\mathrm{d}u.
\end{eqnarray*}
Using this expression, it is easy to show that $|a_n|_\infty\leq
|w'|_\infty/2$. And since $w$ is of class $C^2$, a mean-value theorem
on $w'$ using the above expression shows that $|a_n|_\mathrm{tv}= |a_{n,1}|+
\sum_{i=2}^n|a_{n,i}-a_{n,i-1}|\leq(|w'|_\infty+ |w^{\prime\prime}|_\infty)/2$.
We are then in position to apply Lemma~\ref{techlem1}(1) to obtain
\[
\mathbb{E} \Biggl[ \Biggl(\frac{1}{\sqrt{n}}\sum_{i=1}^na_{n,i}
\bar h(X_i) \Biggr)^2 \Biggr]= \mathrm{O}(1).
\]
By similar arguments as above, and since $u_n=n^{-2}\sum_{i=1}^n \sum_{=1}^n w (\frac{i-j}{n} )$ is a Riemann sum approximation of
$\int_0^1 v(t)\,\mathrm{d}t$, we obtain that $\llvert u_n-\int_0^1v(t)\,\mathrm{d}
t\rrvert =\mathrm{O} (\frac{1}{n} )$. In conclusion,
\begin{equation}
\label{bound-diff1} \mathbb{E} \bigl(\bigl|\sigma_n^2-\check
\sigma_n^2\bigr| \bigr)=\mathrm{O} \biggl(\frac
{1}{n}
\biggr).
\end{equation}
%
\vspace*{2pt}\upqed\end{pf}

\begin{lemma}\label{lemtechlem2thm1}
Assume~\textup{A1}. Suppose that the Markov chain $\{X_n, n\geq0\}$
starts at $X_0=x$ for $x\in\mathsf{X}$ such that (\ref{ergo}) holds. Let
$h\in L_0^2(\pi)$ be such that $\sigma_P^2(h)<\infty$. Then as $n\to
\infty$, $\sigma_n^2-\check\sigma_n^2$ converges in probability to
zero, and $Z^{(n)}\stackrel{\mathsf{w}}{\to}Z$ in $\ell^2(\alpha)$.
\end{lemma}

\begin{pf}
Ergodicity is equivalent to the existence of a successful coupling of
the Markov chain and its stationary copy. More precisely, we can
construct a process $\{(X_n,\tilde{X}_n), n\geq0\}$ such that $\{X_n,
n\geq0\}$ is a Markov chain with initial distribution $\delta_x$ and
transition kernel $P$, $\{\tilde{X}_n, n\geq0\}$ is a Markov chain
with initial distribution $\pi$ and transition kernel $P$, and there
exists a finite (coupling) time $\tau$ such that $X_n=\tilde{X}_n$ for
all $n\geq\tau$. For a proof of this result, see for instance
Lindvall \cite{lindvall92}, Theorem~14.10; see also
Roberts and Rosenthal \cite{robertsetrosenthalsurvey},
Proposition~28. We use a wide ``tilde'' to denote quantities computed
from the stationary chain $\{\tilde{X}_n, n\geq0\}$.

Since $X_n=\tilde{X}_n$ for all $n\geq\tau$, and in view of the
expression of $\sigma_n^2-\check\sigma_n^2$ given in (\ref
{eqdiffsigma}), it is straightforward to check that $\sigma
_n^2-\check
\sigma_n^2- (\widetilde{\sigma_n^2-\check\sigma_n^2} )$
converges to zero in probability. The convergence of $\|
Z^{(n)}-\widetilde{Z^{(n)}}\|_\alpha$ is handled similarly.
\begin{eqnarray*}
&& \bigl\|Z^{(n)}-\widetilde{Z^{(n)}}\bigr\|_\alpha^2\\
 &&\quad=
\sum_{\ell\in\mathsf
{I}}\alpha _\ell \Biggl(
\frac{1}{\sqrt{n}}\sum_{k=1}^n
\Psi_\ell \biggl(\frac
{k}{n} \biggr) \bigl(h(X_k)-h(
\tilde{X}_k) \bigr) \Biggr)^2
\\
&&\quad= \sum_{\ell\in\mathsf{I}}\alpha_\ell \Biggl(
\frac{1}{\sqrt
{n}}\sum_{k=1}^{\tau-1}
\Psi_\ell \biggl(\frac{k}{n} \biggr) \bigl(h(X_k)-h(
\tilde{X}_k) \bigr) \Biggr)^2
\\
&&\quad\leq  \frac{\tau}{n} \biggl(\sup_{t\in[0,1]}\sum
_{\ell\in\mathsf{I}
}\alpha_\ell \bigl|\Psi_\ell(t)\bigr|^2
\biggr) \Biggl(\sum_{k=1}^\tau
\bigl(h(X_k)-h(\tilde{X}_k) \bigr)^2 \Biggr),
\end{eqnarray*}
which converges almost surely to zero, given (\ref{boundPsi}), and
since $\tau$ is finite almost surely.
\end{pf}

\subsection{Proof of Theorem~\texorpdfstring{\protect\ref{thmRWM}}{2.4}}\label{secproof-thmRWM}
Since $u$ is bounded from below, we can choose $a=1 - \inf_{x\in
\mathsf{X}}
u(x)$ such that $V(x)\stackrel{\mathrm{def}}{=}a +u(x)\geq1$. Let
$q_\sigma(x,y)$ be the
density of the proposal $\mathbf{N} (x-\frac{\sigma^2}{2}\rho
(x)\nabla u(x),\sigma^2I_d )$, and define $\mathsf
{R}(x)\stackrel{\mathrm{def}}{=}
\{
y\in\mathbb{R}^p \dvt \alpha(x,y) < 1\}$. We have
\begin{eqnarray}
PV(x)-V(x) &=& \int\alpha(x,y) \bigl(V(y)-V(x)
\bigr)q_\sigma (x,y)\,\mathrm{d} y
\nonumber
\\
\label{PV}
&=& \int_{\mathsf{R}(x)} \bigl[\alpha(x,y)-1 \bigr] \bigl(V(y)-V(x)
\bigr)q_\sigma(x,y)\,\mathrm{d}y
\\
&&{}+ \int \bigl(V(y)-V(x) \bigr)q_\sigma(x,y)\,\mathrm{d}y.\nonumber
\end{eqnarray}
Since $\nabla u$ is Lipschitz, with Lipschitz constant $L$, say, we
have by Taylor expansion
\[
V(y)-V(x)\leq \bigl\anglel \nabla u(x),y-x \bigr\angler _2 +
\frac{L}{2}|y-x|^2.
\]
Integrating both sides, and using the fact that $\rho(x)|\nabla
u(x)|\leq\tau$, we get
\begin{eqnarray}
\int \bigl(V(y)-V(x) \bigr)q_\sigma(x,y)\,\mathrm{d}y
& \leq & -\frac{\sigma^2}{2}\rho(x)\bigl|\nabla u(x)\bigr|^2 +\frac{L}{2}
\biggl(\frac
{\sigma^4}{4}\rho(x)^2\bigl|\nabla u(x)\bigr|^2 +\mathrm{d}
\sigma^2 \biggr)
\nonumber
\\[-8pt]
\label{bound1}
\\[-8pt]
\nonumber
& \leq & -\frac{\sigma^2}{2}\rho(x)\bigl|\nabla u(x)\bigr|^2 +\frac{L}{2}
\biggl(\frac
{\tau^2\sigma^4}{4} +\mathrm{d}\sigma^2 \biggr).
\end{eqnarray}
We also have
\begin{eqnarray*}
&& \frac{\pi(y)}{\pi(x)}\frac{q_\sigma(y,x)}{q_\sigma(x,y)}
\\
&&\quad=\exp \biggl(V(x)-V(y) -\frac{1}{2\sigma^2} \biggl\llvert x-y+\frac{\sigma
^2}{2}
\rho(y)\nabla u(y) \biggr\rrvert ^2 +\frac{1}{2\sigma^2} \biggl\llvert
y-x-\frac
{\sigma^2}{2}\rho(x)\nabla u(x) \biggr\rrvert ^2 \biggr).
\end{eqnarray*}
If $y\in\mathsf{R}(x)$, we necessarily have $\frac{\pi(y)}{\pi
(x)}\frac
{q_\sigma(y,x)}{q_\sigma(x,y)}<1$, which translates to
\[
V(y)-V(x)> -\frac{1}{2\sigma^2} \biggl\llvert x-y+\frac{\sigma^2}{2}\rho (y)
\nabla u(y) \biggr\rrvert ^2 +\frac{1}{2\sigma^2} \biggl\llvert y-x-
\frac{\sigma
^2}{2}\rho(x)\nabla u(x) \biggr\rrvert ^2.
\]
Hence, if $y\in\mathsf{R}(x)$,
\begin{eqnarray*}
&&\bigl[\alpha(x,y)-1 \bigr] \bigl(V(y)-V(x) \bigr)
\\
&&\quad\leq \bigl[\alpha(x,y)-1 \bigr] \biggl(-\frac{1}{2\sigma^2} \biggl\llvert x-y+
\frac{\sigma^2}{2}\rho(y)\nabla u(y) \biggr\rrvert ^2 +
\frac
{1}{2\sigma
^2} \biggl\llvert y-x-\frac{\sigma^2}{2}\rho(x)\nabla u(x) \biggr
\rrvert ^2 \biggr)
\\
&&\quad= \bigl[1-\alpha(x,y) \bigr]\frac{\sigma^2}{8} \biggl(\rho ^2(y)\bigl|
\nabla u(y)\bigr|^2 -\rho^2(x)\bigl|\nabla u(x)\bigr|^2\\
&&\hspace*{74pt}\qquad{}-
\frac{2}{\sigma^2} \bigl\anglel y-x,\rho (x)\nabla u(x) + \rho(y)\nabla u(y) \bigr
\angler \biggr)
\\
&&\quad\leq\frac{\sigma^2}{8} \biggl(\tau^2 +\frac{4\tau}{\sigma
^2}|y-x|
\biggr).
\end{eqnarray*}
Hence,
\begin{equation}
\label{bound2}
\int_{\mathsf{R}(x)} \bigl[\alpha(x,y)-1 \bigr]
\bigl(V(y)-V(x) \bigr)q_\sigma(x,y)\,\mathrm{d}y \leq\frac{\sigma^2\tau^2}{8} +
\frac
{\tau
}{2}\sqrt {\mathrm{d}\sigma^2 +\frac{\sigma^2\tau^2}{2}}.
\end{equation}
We combine (\ref{PV})--(\ref{bound2}) to conclude that
\[
PV(x)-V(x)\leq-\frac{\sigma^2}{2}\rho(x)\bigl|\nabla u(x)\bigr|^2 +K,
\]
where $K=\frac{L}{2} (\frac{\tau^2\sigma^4}{4} +\mathrm{d}\sigma
^2 )+
\frac{\sigma^2\tau^2}{8} +\frac{\tau}{2}\sqrt{\mathrm{d}\sigma^2 +\frac
{\sigma
^2\tau^2}{2}}$. Since $f(x)\stackrel{\mathrm{def}}{=}\frac{\sigma
^2}{2}\rho(x)|\nabla
u(x)|^2$ is continuous and $f(x)\to\infty$, as $\|x\|\to\infty$ by
assumption, the results follow readily.

\subsection{Proof Theorem \texorpdfstring{\protect\ref{thm2}}{2.6}}\label{proofthm2}
We follow Dedecker and Rio \cite{dedeetrio08}, Theorem~2.1. With the geometric ergodicity
assumption, the martingale approximation to $\sum_{i=1}^n h(X_i)$ can
be constructed more explicitly than in Lemmas~\ref{theokv} and~\ref{techlem1}. Define
\[
g(x)=\sum_{j\geq0} P^j\bar h(x),\qquad x\in
\mathsf{X}.
\]
By the geometric ergodicity assumption, $g$ is well-defined and belongs
to $\mathcal{L}_{V^\delta}$. Then we define $D_0=0$, and
$D_k\stackrel{\mathrm{def}}{=}
g(X_k)-Pg(X_{k-1})$, $k\geq1$. It is easy to see that $\{D_k, k\geq
0\}$ is a martingale-difference sequence with respect to the natural
filtration of $\{X_n, n\geq0\}$. Using this martingale, we define
\[
\bar\sigma_n^2\stackrel{\mathrm{def}} {=}\sum
_{\ell\in\mathsf
{I}}\alpha_\ell \Biggl(\frac
{1}{\sqrt
{n}}\sum
_{i=1}^n\Psi_\ell \biggl(
\frac{i-1}{n} \biggr)D_i \Biggr)^2,
\]
and we recall that $\check\sigma_n^2\stackrel{\mathrm{def}}{=}\sum_{\ell\in\mathsf{I}}\alpha
_\ell
 (\frac{1}{\sqrt{n}}\sum_{i=1}^n\Psi_\ell (\frac
{i-1}{n}
)h(X_i) )^2$. Hence,
\[
\sigma_n^2=\sum_{\ell\in\mathsf{I}}
\alpha_\ell \Biggl(\frac
{1}{\sqrt
{n}}\sum_{i=1}^n
\Psi_\ell \biggl(\frac{i-1}{n} \biggr)D_i
\Biggr)^2 + \bigl(\sigma _n^2-\check
\sigma_n^2 \bigr) + \bigl(\check\sigma_n^2-
\bar\sigma _n^2 \bigr).
\]
Although the martingales are constructed differently, the argument in
Lemma~\ref{lemtechlem1thm1} carries through and shows that $\mathbb
{E}
(|\sigma_n^2-\check\sigma_n^2| ) =\mathrm{O}(1/n)$. The proof is
similar to
the proof of Lemma~\ref{lemtechlem1thm1} and is omitted. Also
$\mathbb{E}
(|\check\sigma_n^2-\bar\sigma_n^2| ) =\mathrm{O}(1/\sqrt{n})$. To see this,
use the Cauchy--Schwarz inequalities for sequences in $\ell^2(\alpha)$
and for random variables to write
\begin{eqnarray*}
&& \mathbb{E} \bigl(\bigl|\check\sigma_n^2-\bar
\sigma_n^2\bigr| \bigr)
\\
&&\quad=\mathbb{E} \Biggl[ \Biggl\llvert \sum_{\ell\in\mathsf{I}}
\alpha_\ell \Biggl(\frac
{1}{\sqrt{n}}\sum_{i=1}^n
\Psi_\ell \biggl(\frac{i-1}{n} \biggr) \bigl(h(X_i)-D_i
\bigr) \Biggr) \Biggl(\frac{1}{\sqrt{n}}\sum_{i=1}^n
\Psi_\ell \biggl(\frac
{i-1}{n} \biggr) \bigl(h(X_i)+D_i
\bigr) \Biggr) \Biggr\rrvert \Biggr]
\\
&&\quad\leq\mathbb{E} \Biggl[ \Biggl\{\sum_{\ell\in\mathsf{I}}
\alpha_\ell \Biggl(\frac
{1}{\sqrt
{n}}\sum_{i=1}^n
\Psi_\ell \biggl(\frac{i-1}{n} \biggr) \bigl(h(X_i)-D_i
\bigr) \Biggr)^2 \Biggr\}^{1/2}
\\
&&\hspace*{12pt}\quad\quad{}\times \Biggl\{\sum_{\ell\in\mathsf{I}}\alpha_\ell
\Biggl(\frac
{1}{\sqrt
{n}}\sum_{i=1}^n
\Psi_\ell \biggl(\frac{i-1}{n} \biggr) \bigl(h(X_i)+D_i
\bigr) \Biggr)^2 \Biggr\}^{1/2} \Biggr]
\\
&&\quad\leq \Biggl\{\sum_{\ell\in\mathsf{I}}\alpha_\ell
\mathbb{E} \Biggl[ \Biggl(\frac
{1}{\sqrt
{n}}\sum_{i=1}^n
\Psi_\ell \biggl(\frac{i-1}{n} \biggr) \bigl(h(X_i)-D_i
\bigr) \Biggr)^2 \Biggr] \Biggr\}^{1/2}
\\
&&\quad\quad{}\times \Biggl\{\sum_{\ell\in\mathsf{I}}\alpha_\ell
\mathbb {E} \Biggl[ \Biggl(\frac
{1}{\sqrt
{n}}\sum_{i=1}^n
\Psi_\ell \biggl(\frac{i-1}{n} \biggr) \bigl(h(X_i)+D_i
\bigr) \Biggr)^2 \Biggr] \Biggr\}^{1/2}.
\end{eqnarray*}
By the martingale approximation, we have
\begin{eqnarray*}
\sum_{i=1}^n\Psi_\ell \biggl(
\frac{i-1}{n} \biggr) \bigl(h(X_i)-D_i \bigr) &=&
\Psi_\ell (0)Pg(X_0) - \Psi \biggl(\frac{n-1}{n}
\biggr)Pg(X_n)
\\
&&{}+\sum_{i=2}^n \biggl(\Psi_\ell
\biggl(\frac{i-1}{n} \biggr)-\Psi _\ell \biggl(\frac{i-2}{n}
\biggr) \biggr) Pg(X_{i-1}).
\end{eqnarray*}
The details of these calculations can be found for instance in \cite
{atchadeetcattaneo14}, Proposition A1. It is then easy to show that
\begin{eqnarray*}
&&\sum_{\ell\in\mathsf{I}}\alpha_\ell\mathbb{E} \Biggl[
\Biggl(\sum_{i=1}^n\Psi _\ell
\biggl(\frac{i-1}{n} \biggr) \bigl(h(X_i)-D_i
\bigr) \Biggr)^2 \Biggr]
\\
&&\quad\leq \biggl(6\sup_{0\leq t\leq1}\sum_{\ell\in\mathsf{I}}
\alpha _\ell \bigl|\Psi_\ell (t)\bigr|^2 + 3\sup
_{0\leq t\leq1}\sum_{\ell\in\mathsf{I}}
\alpha_\ell \bigl|\Psi _\ell '(t)\bigr|^2
\biggr)|h|_{V^\delta}^2.
\end{eqnarray*}
For the second term, notice that
\[
\sum_{i=1}^n\Psi_\ell \biggl(
\frac{i-1}{n} \biggr) \bigl(h(X_i)+D_i \bigr) = 2
\sum_{i=1}^n\Psi_\ell \biggl(
\frac{i-1}{n} \biggr)D_i + \sum_{i=1}^n
\Psi _\ell \biggl(\frac{i-1}{n} \biggr) \bigl(h(X_i)-D_i
\bigr).
\]
Hence, with similar calculations, we obtain
\begin{eqnarray*}
&& \sum_{\ell\in\mathsf{I}}\alpha_\ell\mathbb{E} \Biggl[
\Biggl(\sum_{i=1}^n\Psi _\ell
\biggl(\frac{i-1}{n} \biggr) \bigl(h(X_i)+D_i
\bigr) \Biggr)^2 \Biggr] \\
&&\quad\leq   2|h|_{V^\delta
}^2n\sup
_{0\leq t\leq1}\sum_{\ell\in\mathsf{I}}
\alpha_\ell\bigl|\Psi _\ell (t)\bigr|^2
\\
&&\qquad{}+ 6 \biggl(2\sup_{0\leq t\leq1}\sum_{\ell\in\mathsf{I}}
\alpha _\ell\bigl|\Psi _\ell (t)\bigr|^2 + \sup
_{0\leq t\leq1}\sum_{\ell\in\mathsf{I}}
\alpha_\ell \bigl|\Psi _\ell '(t)\bigr|^2
\biggr)|h|_{V^\delta}^2.
\end{eqnarray*}
Given (\ref{boundPsi}), these calculations show that $\mathbb{E}
(|\check
\sigma_n^2-\bar\sigma_n^2| )=\mathrm{O}(1/\sqrt{n})$. We conclude that
\[
\sigma_n^2=\sum_{\ell\in\mathsf{I}}
\alpha_j \Biggl(\frac
{1}{\sqrt
{n}}\sum_{i=1}^n
\Psi_\ell \biggl(\frac{i-1}{n} \biggr)D_i
\Biggr)^2 +\mathrm{O}_p \biggl(\frac
{1}{\sqrt{n}} \biggr),
\]
which implies that
\begin{equation}
\label{step1} \mathsf{d}_1 \bigl(\sigma_n^2,
\chi^2 \bigr)\lesssim\mathsf{d}_1 \bigl(\bar
\sigma_n^2,\chi^2 \bigr) + \frac{1}{\sqrt{n}}.
\end{equation}
Therefore, we only need to focus on the term $\mathsf{d}_1 (\bar
\sigma_n^2,\chi^2 )$.

On the Euclidean space $\mathbb{R}^{\mathsf{I}}$, we define the norms
$\|x\|
_\alpha
^2=\sum_{i\in\mathsf{I}}\alpha_ix_i^2$, $\|x\|^2=\sum_{i\in
\mathsf{I}}x_i^2$ and the
inner-products $\anglel x,y\angler _\alpha=\sum_{i\in\mathsf
{I}}\alpha_ix_iy_i$,
and $\anglel x,y\angler =\sum_{i\in\mathsf{I}}x_iy_i$. For a sequence
$(a_1,a_2,\ldots)$, we use
the notation $a_{i:k}=(a_i,\ldots,a_k)$ (and $a_{i:k}$ is the empty set
if $i>k$).
We introduce new random variables $\{Z_{i,j},  i\in\mathsf{I},
1\leq j\leq
n\}
$ which are i.i.d. $\mathbf{N}(0,1)$, and set $S_{\ell:k}\stackrel
{\mathrm{def}}{=}
(\sum_{j=\ell}^k Z_{1j},\ldots,\sum_{j=\ell}^k Z_{\mathsf
{I}j} )^\mathsf{T}
\in
\mathbb{R}^{\mathsf{I}}$, so that
\[
\chi^2\stackrel{\mathrm{dist.}} {=} \sum_{i\in\mathsf{I}}
\alpha_i \Biggl(\frac
{1}{\sqrt
{n}}\sum_{j=1}^n
Z_{i,j} \Biggr)^2=\biggl\|\frac{1}{\sqrt{n}}S_{1:n}\biggr\|
_\alpha^2.
\]
For $1\leq\ell\leq k\leq n$, and omitting the dependence on $n$, we
set ${\mathbf B}_{\ell:k}$ as the $\mathbb{R}^{\mathsf{I}\times(k-\ell
+1)}$ matrix
\[
{\mathbf B}_{\ell:k}(i,j)=\Psi_i \biggl(\frac{j}{n}
\biggr),\qquad i\in \mathsf{I}, \ell \leq j\leq k.
\]
By the Mercer's expansion for $\phi$, we have
\[
\bar\sigma_n^2=\sum_{i\in\mathsf{I}}
\alpha_i \Biggl(\frac
{1}{\sqrt
{n}}\sum_{k=1}^n
\Psi_i \biggl(\frac{k}{n} \biggr)D_k
\Biggr)^2=\biggl\|\frac
{1}{\sqrt
{n}}{\mathbf B}_{1:n}D_{1:n}
\biggr\|_\alpha^2.
\]
For $f\in\mathsf{Lip}_1(\mathbb{R})$, we introduce the function
$f_\alpha \dvtx  \mathbb{R}
^{|\mathsf{I}|}\to\mathbb{R}$, defined as $f_\alpha(x)=f (\|x\|
_\alpha
^2
)$. As a matter of telescoping the sums, we have
\begin{eqnarray*}
&& \mathbb{E} \bigl[f \bigl(\bar\sigma^2_n \bigr)-f \bigl(
\chi^2 \bigr) \bigr]\\
&&\quad= \mathbb {E} \biggl[f_\alpha \biggl(
\frac
{1}{\sqrt{n}}{\mathbf B}_{1:n}D_{1:n}
\biggr)-f_\alpha \biggl(\frac
{1}{\sqrt
{n}}S_{1:n} \biggr)
\biggr]
\\
&&\quad= \sum_{\ell=1}^n \mathbb{E}
\biggl[f_\alpha \biggl(\frac{1}{\sqrt
{n}}{\mathbf B}_{1:\ell}D_{1:\ell}+
\frac{1}{\sqrt{n}}S_{\ell+1:n} \biggr)-f_\alpha \biggl(
\frac{1}{\sqrt{n}}{\mathbf B}_{1:\ell-1}D_{1:\ell-1}+
\frac
{1}{\sqrt
{n}}S_{\ell:n} \biggr) \biggr]
\\
&&\quad= \sum_{\ell=1}^n\mathbb{E}
\biggl[f_{\alpha,n,\ell+1} \biggl(\frac
{1}{\sqrt
{n}}{\mathbf B}_{1:\ell-1}D_{1:\ell-1}+
\frac{1}{\sqrt{n}}{\mathbf B}_{\ell
:\ell
}D_\ell
\biggr)\\
&&\hspace*{27pt}\qquad{}-f_{\alpha,n,\ell+1} \biggl(\frac{1}{\sqrt{n}}{\mathbf B}_{1:\ell-1}D_{1:\ell-1}+
\frac{1}{\sqrt{n}}S_{\ell:\ell} \biggr) \biggr],
\end{eqnarray*}
where we define
\[
f_{\alpha,n,\ell}(x)\stackrel{\mathrm{def}} {=}\mathbb{E} \biggl[f_\alpha
\biggl(x+\frac
{1}{\sqrt
{n}}S_{\ell:n} \biggr) \biggr]\quad \mbox{and set }
f_{\alpha
,n,n+1}(x)=f_\alpha(x).
\]
%

First, we claim that $f_{\alpha,n,\ell}$ is differentiable everywhere
on $\mathbb{R}^{\mathsf{I}}$. To prove this, it suffices to obtain
the almost
everywhere differentiability of $z\in\mathbb{R}^{\mathsf{I}}\mapsto
f_\alpha
(x+z )$ for any $x\in\mathbb{R}^\mathsf{I}$. By Rademacher's
theorem, $f$
as a
Lipschitz function is differentiable almost everywhere on $\mathbb
{R}$. If
$E$ is the set of points where $f$ is not differentiable, we conclude
that $f_\alpha$ is differentiable at all points $z\notin\{z\in
\mathbb{R} ^{\mathsf{I}
}\dvt  \|x+z\|_\alpha^2\in E\}$. Now by\vspace*{1.5pt} Ponomar{\"e}v \cite{ponomarev87}, Theorem~2, the
Lebesgue measure of the set $\{z\in\mathbb{R}^{\mathsf{I}} \dvt \|x+z\|
_\alpha^2\in
E\}
$ is zero, which proves the claim.

As a result, the function $x\mapsto f_{\alpha,n,\ell}(x)$ is
differentiable with derivative
\[
\nabla f_{\alpha,n,\ell}(x)\cdot h=2\mathbb{E} \biggl[f'_\alpha
\biggl(x+\frac
{1}{\sqrt{n}}S_{\ell:n} \biggr) \biggl\anglel x+
\frac{1}{\sqrt{n}}S_{\ell
:n},h \biggr\angler _\alpha \biggr].
\]
By writing this expectation wrt the distribution of $x+\frac{1}{\sqrt
{n}}S_{\ell:n}$, we get
\[
\nabla f_{\alpha,n,\ell}(x)\cdot h=2\int f'_\alpha(z)
\anglel z,h\angler _\alpha \exp \biggl(-\frac{n}{2(n-\ell+1)} \bigl(\|x
\|^2-2\anglel x,z\angler \bigr) \biggr)\mu_{n,\ell}(
\mathrm{d}z),
\]
where $\mu_{n,\ell}$ is the distribution of $\frac{1}{\sqrt
{n}}S_{\ell
:n}$. This implies that $f_{\alpha,n,\ell}$ is infinitely
differentiable with second derivatives given by
\begin{eqnarray*}
&& \nabla^{(2)} f_{\alpha,n,\ell}(x)\cdot(h_1,h_2)
\\
&&\quad=-2 \biggl(\frac{n}{n-\ell+1} \biggr)\\
&&\qquad{}\times\int f'_\alpha(z)
\anglel z,h_1\angler _\alpha \anglel x-z,h_2
\angler \exp \biggl(-\frac{n}{2(n-\ell+1)} \bigl(\|x\| ^2-2\anglel x,z
\angler \bigr) \biggr)\mu_{n,\ell}(\mathrm{d}z)
\\
&&\quad=2\mathbb{E} \biggl[f'_\alpha \biggl(x+
\frac{1}{\sqrt{n}}S_{\ell
:n} \biggr) \biggl\anglel x\sqrt{
\frac{n}{n-\ell+1}}+\frac{S_{\ell:n}}{\sqrt
{n-\ell +1}},h_1 \biggr\angler
_\alpha \biggl\anglel \frac{S_{\ell:n}}{\sqrt
{n-\ell +1}},h_2 \biggr\angler
\biggr],
\end{eqnarray*}
which implies after some easy calculations that for $h\in\mathbb
{R}^\mathsf{I}$,
\begin{equation}
\label{boundD2}
\bigl\llvert \nabla^{(2)} f_{\alpha,n,\ell}(x)\cdot(h,h)
\bigr\rrvert \lesssim\|h\|^2 \biggl(1+\sqrt{\frac{n}{n-\ell+1}}\|x
\|_\alpha \biggr).
\end{equation}
Similarly
for $h\in\mathbb{R}^\mathsf{I}$,
\begin{equation}
\label{boundD3}
\bigl\llvert \nabla^{(3)} f_{\alpha,n,\ell}(x)
\cdot(h,h,h) \bigr\rrvert \lesssim\sqrt{\frac{n}{n-\ell+1}}\|h\|^3
\biggl(1+\sqrt{\frac
{n}{n-\ell
+1}}\|x\|_\alpha \biggr).
\end{equation}

Now, by Taylor expansion we have
\begin{eqnarray*}
&& f_{\alpha,n,\ell+1} \biggl(\frac{1}{\sqrt{n}}{\mathbf B}_{1:\ell
-1}D_{1:\ell
-1}+
\frac{1}{\sqrt{n}}{\mathbf B}_{\ell:\ell}D_\ell
\biggr)-f_{\alpha
,n,\ell
+1} \biggl(\frac{1}{\sqrt{n}}{\mathbf B}_{1:\ell-1}D_{1:\ell-1}+
\frac
{1}{\sqrt
{n}}S_{\ell:\ell} \biggr)
\\
&&\quad=\frac{1}{\sqrt{n}}\nabla f_{\alpha,n,\ell+1} \biggl(\frac{1}{\sqrt
{n}}{\mathbf
B}_{1:\ell-1}D_{1:\ell-1} \biggr)\cdot ({\mathbf B}_{\ell:\ell}D_\ell
-S_{\ell
:\ell} )
\\
&&\quad\quad{}+\frac{1}{2n}\nabla^{(2)} f_{\alpha,n,\ell+1} \biggl(
\frac{1}{\sqrt
{n}}{\mathbf B}_{1:\ell-1}D_{1:\ell-1} \biggr)\cdot \bigl[
({\mathbf B}_{\ell:\ell
}D_\ell ,{\mathbf B}_{\ell:\ell}D_\ell
)- (S_{\ell:\ell},S_{\ell
:\ell
} ) \bigr] + \varrho_{n,\ell}^{(3)},
\end{eqnarray*}
where, using (\ref{boundD3}),
\[
\bigl\llvert \varrho_{n,\ell}^{(3)} \bigr\rrvert \lesssim\sqrt{
\frac{n}{n-\ell
+1}}n^{-3/2} \biggl(1+\sqrt{\frac{\ell-1}{n-\ell+1}} \biggl
\llVert \frac
{{\mathbf
B}_{1:\ell-1}D_{1:\ell-1}}{\sqrt{\ell-1}} \biggr\rrVert _\alpha \biggr) \bigl(\| {\mathbf
B}_{\ell:\ell}D_\ell\|_\alpha^3+
\|S_{\ell:\ell}\|_\alpha ^3 \bigr).
\]
It follows that
\begin{equation}
\label{boundrho3} \sum_{\ell=1}^{n-1}\mathbb{E}
\bigl( \bigl\llvert \varrho_{n,\ell
}^{(3)} \bigr\rrvert \bigr)
\lesssim n^{-1}\sum_{\ell=1}^{n}
\frac{1}{\sqrt{\ell}} + n^{-1/2}\sum_{\ell=1}^n
\frac{1}{\ell}\lesssim n^{-1/2}\log(n).
\end{equation}
By first conditioning on $\mathcal{F}_{\ell-1}$, we have
\[
\mathbb{E} \biggl[\nabla f_{\alpha,n,\ell+1} \biggl(\frac{1}{\sqrt{n}}{\mathbf
B}_{1:\ell
-1}D_{1:\ell-1} \biggr)\cdot ({\mathbf B}_{\ell:\ell}D_\ell-S_{\ell
:\ell}
) \biggr]=0.
\]

Writing $K_{n,\ell}\stackrel{\mathrm{def}}{=}\frac{1}{2}\nabla
^{(2)} f_{\alpha,n,\ell
}(\frac{1}{\sqrt{n}}{\mathbf B}_{1:\ell-1}D_{1:\ell-1})$, we have
\begin{eqnarray*}
&& \nabla^{(2)} f_{\alpha,n,\ell+1} \biggl(\frac{1}{\sqrt{n}}{\mathbf
B}_{1:\ell
-1}D_{1:\ell-1} \biggr)\cdot \bigl[ ({\mathbf
B}_{\ell:\ell}D_\ell,{\mathbf B}_{\ell
:\ell}D_\ell )-
(S_{\ell:\ell},S_{\ell:\ell} ) \bigr]
\\
&&\quad=D_\ell^2\sum_{i,j}
\Psi_i \biggl(\frac{\ell}{n} \biggr)\Psi_j \biggl(
\frac
{\ell}{n} \biggr)K_{n,\ell}(i,j)- \sum
_{i,j}\Psi_i \biggl(\frac
{\ell
}{n} \biggr)
\Psi_j \biggl(\frac{\ell}{n} \biggr)K_{n,\ell
}(i,j)Z_{i,\ell
}Z_{j\ell}.
\end{eqnarray*}
Therefore,
\begin{eqnarray*}
&&\mathbb{E} \biggl(\nabla^{(2)} f_{\alpha,n,\ell+1} \biggl(
\frac{1}{\sqrt
{n}}{\mathbf B}_{1:\ell-1}D_{1:\ell-1} \biggr)\cdot \bigl[
({\mathbf B}_{\ell:\ell
}D_\ell ,{\mathbf B}_{\ell:\ell}D_\ell
)- (S_{\ell:\ell},S_{\ell
:\ell
} ) \bigr]\vert\mathcal{F}_{\ell-1}
\biggr)
\\
&&\quad =\sum_{i,j}\Psi_i \biggl(
\frac{\ell}{n} \biggr)\Psi_j \biggl(\frac
{\ell
}{n}
\biggr)K_{n,\ell+1}(i,j) \bigl[\mathbb{E} \bigl(D_\ell^2
\vert \mathcal{F} _{\ell
-1} \bigr)-\delta_{ij} \bigr],
\end{eqnarray*}
where $\delta_{ij}=1$ if $i=j$ and zero otherwise. We claim that the
proof will be finished if we show that for all $i,j\in\mathsf{I}$,
and $1\leq
\ell\leq n$,
\begin{equation}
\label{claim2}
\mathbb{E}^{1/2} \bigl[ \bigl(K_{n,\ell}(i,j)-K_{n,\ell
+1}(i,j)
\bigr)^2 \bigr]\lesssim\frac{\sqrt{n}}{n-\ell+1}.
\end{equation}
To prove this claim, it suffice to use (\ref{claim2}) to show that
$\llvert n^{-1}\sum_{\ell=1}^n\Psi_i (\frac{\ell}{n}
)\Psi_j
(\frac{\ell}{n} )\mathbb{E} (K_{n,\ell+1}(i,j)
)\rrvert \lesssim
n^{-1/2}\log(n)$ for $i\neq j$, and $\llvert n^{-1}\sum_{\ell
=1}^n\Psi
_i (\frac{\ell}{n} )\Psi_j (\frac{\ell}{n}
)\mathbb{E}
(K_{n,\ell+1}(i,j) [\mathbb{E} (D_\ell^2\vert\mathcal
{F}_{\ell-1}
)-1
] )\rrvert \lesssim n^{-1/2}\log(n)$ for all $i,j\in\mathsf{I}$.
To show this, write
\begin{eqnarray*}
&& \frac{1}{n}\sum_{\ell=1}^{n-1}
\Psi_i \biggl(\frac{\ell}{n} \biggr)\Psi _j \biggl(
\frac{\ell}{n} \biggr)\mathbb{E} \bigl(K_{n,\ell
+1}(i,j) \bigr)\\
&&\quad= \Biggl\{
\frac{1}{n}\sum_{\ell=1}^{n-1}
\Psi_i \biggl(\frac{\ell}{n} \biggr)\Psi _j \biggl(
\frac{\ell}{n} \biggr) \Biggr\}\mathbb{E} \bigl(K_{n,n}(i,j) \bigr)
\\
&&\qquad{}+\frac{1}{n}\sum_{\ell=1}^{n-1} \Biggl[
\frac{1}{n}\sum_{k=1}^{\ell
-1}\Psi
_i \biggl(\frac{\ell}{n} \biggr)\Psi_j \biggl(
\frac{\ell}{n} \biggr) \Biggr] \bigl[\mathbb{E} \bigl(K_{n,\ell}(i,j)-K_{n,\ell+1}(i,j)
\bigr) \bigr].
\end{eqnarray*}
By the convergence of Riemann sums, $\llvert \frac{1}{n}\sum_{\ell
=1}^{n-1}\Psi_i (\frac{\ell}{n} )\Psi_j (\frac
{\ell
}{n} )\rrvert \lesssim n^{-1}$. Combined with (\ref{boundD2}) and (\ref{claim2}), this implies that
\[
\Biggl\llvert \frac{1}{n}\sum_{\ell=1}^n
\Psi_i \biggl(\frac{\ell
}{n} \biggr)\Psi _j \biggl(
\frac{\ell}{n} \biggr)\mathbb{E} \bigl(K_{n,\ell
+1}(i,j) \bigr) \Biggr
\rrvert \leq\frac{1}{n} \Biggl(\sqrt{n}+\sqrt{n}\sum
_{k=1}^n\frac
{1}{k} \Biggr)\lesssim
\frac{\log(n)}{\sqrt{n}}.
\]

For the second term, notice from the definition of $D_\ell$ at the
beginning of the proof that $\mathbb{E} (D_\ell^2\vert\mathcal
{F}_{\ell
-1}
)-1=G(X_{\ell-1})-\pi(G)$, where $G(x)=Pg^2(x)-(Pg(x))^2$. Since
$h\in\mathcal{L}
_{V^\delta}$ for $\delta<1/4$, $G\in\mathcal{L}_{V^{2\delta}}$, and
$2\delta
<1/2$. Therefore, by geometric ergodicity, the solution of the Poisson
equation for $G$ defined as $U(x)=\sum_{j\geq0}P^j(G(x)-\pi(G))$ is
well-defined, $U\in\mathcal{L}_{V^{2\delta}}$, and we have almost surely
\[
U(X_{\ell-1})-PU(X_{\ell-1})=\mathbb{E} \bigl(D_\ell^2
\vert \mathcal{F}_{\ell
-1} \bigr)-1.
\]
Notice that, since $2\delta<1/2$, for any $p\geq2$ such that
$2p\delta
\leq1$, the geometric ergodicity assumption (G) implies that $\sup_{k\geq1} \mathbb{E} (|U(X_k)|^p )<\infty$. Now we use
the usual
martingale approximation trick (see, e.g.,
Atchad{\'e} and Cattaneo
\cite{atchadeetcattaneo14}, Proposition~A1) to write
\begin{eqnarray*}
&&\frac{1}{n}\sum_{\ell=1}^{n-1}
\Psi_i \biggl(\frac{\ell}{n} \biggr)\Psi _j \biggl(
\frac{\ell}{n} \biggr)\mathbb{E} \bigl(K_{n,\ell
+1}(i,j) \bigl[\mathbb{E}
\bigl(D_\ell^2\vert\mathcal{F}_{\ell-1} \bigr)-1
\bigr] \bigr)\\
&&\quad=\frac
{1}{n}\Psi_i \biggl(\frac{1}{n}
\biggr)\Psi_j \biggl(\frac{1}{n} \biggr)\mathbb{E}
\bigl(K_{n,2}(i,j)U(X_0) \bigr)
\\
&&\quad\quad{}-\frac{1}{n}\Psi_i \biggl(1-\frac{1}{n} \biggr)
\Psi_j \biggl(1-\frac
{1}{n} \biggr)\mathbb{E}
\bigl(K_{n,n}(i,j)U(X_{n-1}) \bigr)
\\
&&\quad\quad{}+\frac{1}{n}\sum_{\ell=1}^{n-1}
\mathbb{E} \biggl[ \biggl\{\Psi _i \biggl(\frac{\ell
}{n} \biggr)
\Psi_j \biggl(\frac{\ell}{n} \biggr)K_{n,\ell
+1}(i,j)\\
&&\hspace*{13pt}\qquad\qquad\qquad{}-\Psi
_i \biggl(\frac{\ell-1}{n} \biggr)\Psi_j \biggl(
\frac{\ell
-1}{n} \biggr)K_{n,\ell}(i,j) \biggr\}U(X_{\ell-1})
\biggr].
\end{eqnarray*}
We now use the fact that $\Psi_i\Psi_j$ is of class $\mathsf{C}^1$ (see
Theorem~\ref{theomercer}(ii)), (\ref{boundD2}), and (\ref{claim2}) to
conclude that
\begin{eqnarray*}
&& \Biggl\llvert \frac{1}{n}\sum_{\ell=1}^{n-1}
\Psi_i \biggl(\frac{\ell
}{n} \biggr)\Psi_j \biggl(
\frac{\ell}{n} \biggr)\mathbb{E} \bigl(K_{n,\ell
+1}(i,j) \bigl[\mathbb{E}
\bigl(D_\ell^2\vert\mathcal{F}_{\ell-1} \bigr)-1
\bigr] \bigr) \Biggr\rrvert
\\
&&\quad\lesssim\frac{1}{\sqrt{n}} + \frac{1}{n}\sum_{\ell
=1}^{n-1}
\mathbb{E} ^{1/2} \bigl( \bigl\llvert K_{n,\ell+1}(i,j)-K_{n,\ell+2}(i,j)
\bigr\rrvert ^2 \bigr)\lesssim\frac{\log(n)}{\sqrt{n}}.
\end{eqnarray*}
This proves the claim. It remains to establish (\ref{claim2}). Write
$\mathbb{E}_\ell$ to denote the expectation operator wrt
$n^{-1/2}{\mathbf
S}_{\ell
:n}$. We then have for any $h_1,h_2\in\mathbb{R}^{\mathsf{I}}$,
\begin{eqnarray*}
&& 2K_{n,\ell}\cdot(h_1,h_2)\\
&&\quad=\nabla^{(2)}f_{\alpha,n,\ell}
\biggl(\frac
{1}{\sqrt{n}}{\mathbf B}_{1:\ell-1}D_{1:\ell-1} \biggr)
\cdot(h_1,h_2)
\\
&&\quad=2 \biggl(\frac{n}{n-\ell+1} \biggr)\\
&&\qquad{}\times\mathbb{E}_\ell
\biggl[f'_\alpha \biggl(\frac
{1}{\sqrt{n}}{\mathbf
B}_{1:\ell-1}D_{1:\ell-1}+\frac{S_{\ell
:n}}{\sqrt
{n}} \biggr) \biggl\anglel
\frac{1}{\sqrt{n}}{\mathbf B}_{1:\ell-1}D_{1:\ell
-1}+
\frac {S_{\ell:n}}{\sqrt{n}},h_1 \biggr\angler _\alpha \biggl\anglel
\frac
{S_{\ell:n}}{\sqrt {n}},h_2 \biggr\angler \biggr]
\\
&&\quad = \biggl(\frac{n-\ell}{n-\ell+1} \biggr)\nabla^{(2)}f_{\alpha
,n,\ell
+1}
\biggl(\frac{1}{\sqrt{n}}{\mathbf B}_{1:\ell-1}D_{1:\ell-1}+
\frac
{S_\ell
}{\sqrt{n}} \biggr)\cdot(h_1,h_2) \\
&&\qquad{}+ \biggl(
\frac{n}{n-\ell+1} \biggr)\mathrm{O} \biggl(\frac{1}{\sqrt{n}} \biggr).
\end{eqnarray*}
Therefore,
\begin{eqnarray*}
&& 2 (K_{n,\ell}-K_{n,\ell+1} )\cdot(h_1,h_2)
\\
&&\quad=\nabla^{(2)}f_{\alpha,n,\ell+1} \biggl(\frac{1}{\sqrt{n}}{\mathbf
B}_{1:\ell
-1}D_{1:\ell-1}+\frac{S_\ell}{\sqrt{n}} \biggr)
\cdot(h_1,h_2) \\
&&\qquad{}-\nabla^{(2)}f_{\alpha,n,\ell+1}
\biggl(\frac{1}{\sqrt{n}}{\mathbf B}_{1:\ell
-1}D_{1:\ell-1}+
\frac{{\mathbf B}_\ell D_\ell}{\sqrt{n}} \biggr)\cdot (h_1,h_2)
\\
&&\qquad{}-\frac{1}{n-\ell+1}\nabla^{(2)}f_{\alpha,n,\ell+1} \biggl(
\frac
{1}{\sqrt
{n}}{\mathbf B}_{1:\ell-1}D_{1:\ell-1}+\frac{S_\ell}{\sqrt{n}}
\biggr)\cdot(h_1,h_2) + \biggl(\frac{n}{n-\ell+1} \biggr)
\mathrm{O} \biggl(\frac{1}{\sqrt{n}} \biggr)
\\
&&\quad=\nabla^{(3)}f_{\alpha,n,\ell+1} \biggl(\frac{1}{\sqrt{n}}{\mathbf
B}_{1:\ell
-1}D_{1:\ell-1}+t\frac{S_\ell}{\sqrt{n}}+(1-t)\frac{{\mathbf B}_\ell
D_\ell
}{\sqrt{n}}
\biggr)\cdot \biggl(h_1,h_2,\frac{S_\ell}{\sqrt
{n}}-
\frac{{\mathbf
B}_\ell D_\ell}{\sqrt{n}} \biggr)
\\
&&\qquad{}-\frac{1}{n-\ell+1}\nabla^{(2)}f_{\alpha,n,\ell+1} \biggl(
\frac
{1}{\sqrt
{n}}{\mathbf B}_{1:\ell-1}D_{1:\ell-1}+\frac{S_\ell}{\sqrt{n}}
\biggr)\cdot(h_1,h_2) + \biggl(\frac{n}{n-\ell+1} \biggr)
\mathrm{O} \biggl(\frac{1}{\sqrt{n}} \biggr),
\end{eqnarray*}
for some $t\in(0,1)$. Using (\ref{boundD2}) and (\ref{boundD3}),
(\ref
{claim2}) follows from the above.

\begin{appendix}
\section*{Appendix: Mercer's theorem}\label{secappendix}
We recall Mercer's theorem below. Part~(i) is the standard Mercer's
theorem, and part (ii) is a special case of a result due to T. Kadota
(Kadota \cite{kadota67}).
%

\setcounter{theo}{0}
\begin{theo}[(Mercer's theorem)]\label{theomercer}
\textup{(i)} Let $k\dvtx  [0,1]\times[0,1]\to\mathbb{R}$ be a continuous
positive semidefinite kernel. Then there exist nonnegative numbers $\{
\lambda_j, j\geq0\}$, and orthonormal functions $\{\phi_j, j\geq
0\}
$, $\phi_j\in L^2([0,1])$, such that $\int_0^1k(x,y)\phi
_j(y)\,\mathrm{d}y=\lambda
_j\phi_j(x)$ for all $x\in[0,1] $, $j\geq0$, and
\[
\lim_{n\to\infty}\sup_{x,y\in[0,1] } \Biggl\llvert k(x,y)-
\sum_{j=0}^n\lambda _j
\phi_j(x)\phi_j(y) \Biggr\rrvert =0.
\]
Furthermore, if $\lambda_j\neq0$, $\phi_j$ is continuous.

\textup{(ii)} Let $k$ as above. If in addition $k$ is of class $\mathsf
{C}^2$ on $[0,1]\times[0,1]$, then for $\lambda_j\neq0$, $\phi_j$ is
of class $\mathsf{C}^1$ on $[0,1]$ and
\[
\lim_{n\to\infty}\sup_{x,y\in[0,1]} \Biggl\llvert
\frac{\partial
^2}{\partial
x\partial y}k(x,y)-\sum_{j=0}^n
\lambda_j\phi_j'(x)\phi_j'(y)
\Biggr\rrvert =0.
\]
\end{theo}

By setting $x=y$, in both expansions, it follows that
\setcounter{equation}{0}
\begin{equation}
\label{boundnormphi} \sup_{0\leq x\leq1} \sum_{j\geq0}
\lambda_j\bigl|\phi_j(x)\bigr|^2\leq\sup
_{0\leq x\leq1}k(x,x)<\infty
\end{equation}
and
\begin{equation}
\label{boundnormphip}
\sup_{0\leq x\leq1} \sum_{j\geq0}
\lambda_j\bigl|\phi_j'(x)\bigr|^2\leq\sup
_{0\leq x\leq1} \biggl\llvert \frac{\partial^2}{\partial u\partial
v}k(u,v)\vert
_{u=x,v=x} \biggr\rrvert <\infty.
\end{equation}
\end{appendix}

\section*{Acknowledgments}
I am grateful to Matias Cattaneo and Shane
Henderson for very helpful discussions.
This work is partly supported by NSF Grants DMS-09-06631,
NSF-SES 1229261 and DMS-12-28164.


\printhistory
\end{document}